\documentclass[]{interact}

\usepackage{epstopdf}
\usepackage[caption=false]{subfig}
\usepackage[numbers,sort&compress]{natbib}
\bibpunct[, ]{[}{]}{,}{n}{,}{,}
\renewcommand\bibfont{\fontsize{10}{12}\selectfont}

\theoremstyle{plain}

\theoremstyle{definition}

\theoremstyle{remark}
\newtheorem{remark}{Remark}

\usepackage{amsmath}
\usepackage{caption}
\usepackage{amssymb,amsmath,amsthm,mathrsfs}
\usepackage[mathscr]{euscript}
\usepackage{graphicx}
\articletype{}

\title{Hurricane Simulation and Nonstationary Extremal Analysis for a Changing Climate}

\author{
	\name{Meagan Carney\textsuperscript{a} Holger Kantz\textsuperscript{b} Matthew Nicol\textsuperscript{c}\thanks{\textsuperscript{a}m.carney@uq.edu.au; \textsuperscript{b}kantz@pks.mpg.de; \textsuperscript{c}nicol@math.uh.edu}}
	\affil{\textsuperscript{a}The University of Queensland, Brisbane, AUS; \textsuperscript{b}Max Planck Institute for the Physics of Complex Systems, Dresden, DE; \textsuperscript{c}University of Houston, Houston, USA}
}
\begin{document}

\maketitle

\begin{abstract}

Particularly important to hurricane risk assessment for coastal regions is finding accurate approximations of return probabilities of maximum windspeeds. Since extremes in maximum windspeed have a direct relationship to minimums in the central pressure, accurate windspeed return estimates rely heavily on proper modeling of the central pressure minima. Using the HURDAT2 database, we show that the central pressure minima of hurricane events can be appropriately modeled by a nonstationary extreme value distribution. We also provide and validate a Poisson distribution with a nonstationary rate parameter to model returns of hurricane events. Using our nonstationary models and numerical simulation techniques from established literature, we perform a simulation study to model returns of maximum windspeeds of hurricane events along the North Atlantic Coast. We show that our revised model agrees with current data and results in an expectation of higher maximum windspeeds for all regions along the coast with the highest maximum windspeeds occurring in the northern part of the coast.

\end{abstract}


\section{Introduction}

Hurricanes and tropical storms bring massive societal impacts and cause economic instabilities. Known for their high windspeeds and downpours, these storms are often accompanied by flooding, wind damage, and travel hazards that lead to large-scale evacuations and a national emergency response. Talk of climate change in recent years and more frequent observations of extreme weather events has inspired research into techniques that provide more accurate estimates of returns and return times of extremes \cite{CAN,CK,exbook}.

Particularly important in hurricane risk assessment for coastal regions is finding accurate approximations of the return probabilities of maximum windspeeds. There have been several studies surrounding maximum windspeed return estimates of hurricanes occurring along the North Atlantic Coast \cite{B,H,CC,SHW,VT}. Many of these studies use the retired HURDAT database which has since been discounted as an unreliable source for future prediction modeling. In 1999, Casson and Coles purposed a hurricane model that allows for approximations of maximum windspeed returns using the tracks and central pressure minima \cite{CC}. The advantage of a model over raw data analysis is that a large number of hurricanes can be simulated to provide more accurate estimates of the tail probabilities and longer year returns of such rare events. The simulation results of this model are in good agreement with the other models and analyses of that decade. However, our findings suggest that this model does not hold up in accuracy when fitted to the updated HURDAT2 database. These inaccuracies can be almost entirely attributed to systematic trends in the observed central pressure and frequency of hurricane events over recent years.

Although there are many factors in a hurricane event that affect the maximum windspeed, we find that the most influential for risk assessment are the central pressure minima and translational velocity of a hurricane at the time of impact with the coast.  Since the central pressure minima have a direct relationship to the windspeed maxima, a better estimate of their probability distribution can provide more accurate returns of extreme highs of windspeed maxima along the coast. Models of an extreme (e.g. minima or maxima) most often take the form of an extreme value distribution \cite{C,exbook}. These distributions have been studied extensively; however, revisions for more complex data analysis settings are often required. 

Following the work in \cite{CC}, we show that we can still reliably model the central pressure minima of a hurricane event using the generalized extreme value distribution (GEV); however, a previously unobserved time dependent trend in the central pressure minima requires adaptations in both the model and methodology. We also provide evidence for a Poisson distribution with a time-dependent rate parameter to model the number of yearly hurricane events that continues into the modern era (post 1965) which previous literature has assumed to be stationary. 

Our revised model results in two major differences in the simulation of coastal risk analysis of hurricane events: 1.) higher maximum windspeeds are expected for all regions along the North Atlantic Coast, including the Gulf Coast and 2.) the highest maximum windspeeds are expected to occur in the northern part of the coast. Higher maximum windspeeds are likely due to a combination of the central pressure minima time dependence and increase in the number of observed hurricane events incorporated into the model. The second observation is arguably more surprising since the number of hurricane events hitting the coast in the north is much lower than regions near the Gulf of Mexico. An increase in the translational velocity as hurricanes travel northward explains this effect.

\section{Methodology}
\subsection{The Wind Field Model, Maximum Windspeeds, and Minimum Pressure}
We describe the Wind Field Model introduced in \cite{NOAA} and the relationship between maximum windspeeds and minimum central pressure.

Given the center location $(\phi_t,\psi_t)$ in the usual geographic coordinates (degrees) and central pressure $p_t$ in hPa of a hurricane measured at the eye at time $t$, the Wind Field Model \cite{NOAA} allows us to model the stochastic process of maximum windspeeds of a hurricane as a sequence of random variables sampled at any given time $t$ by, 
\begin{equation}\label{wfm}
V(R_{\max},\phi_t, p_t, u_t) = 0.865\bigg(K\sqrt{\Delta p_t}-\frac{R_{\max}(\phi_t)f}{2}\bigg)+0.5u_t
\end{equation}
where $K$ is a constant given in $m (s~\text{hPa}^{1/2})^{-1}$, $f = \omega \sin \phi_t$ is the Coriolis parameter $\omega = 7.2982\times 10^{-4}~s^{-1}$, $\Delta p_t = 0.75(1013-p_t)$ is the pressure differential, $R_{\max}(\phi_t)$ is the radius to maximum windspeeds in $m$ sampled from the distribution in the Appendix \ref{rmax}, and $u_t$ is the translational velocity in $m/s$ (meters per second) at time $t$. The translational velocity $u_t$ at a time $t$ is estimated as the change in the distance of the center of the hurricane $\sqrt{(\phi_t-\phi_{t-1})^2+(\psi_t-\psi_{t-1})^2}$ over the change in time $t-1$ to time $t$. For more information on how the variables in the Wind Field Model are related see Table \ref{t:wfm}.

\begin{table}[ht]
	\centering
	\begin{tabular}{|c|c|c|c|c|}
		\hline
		variable: & $V$ & $R_{\max}$ & $p_t$ & $\phi_t $\\
		\hline
		depends on: & $R_{\max}$, $\phi_t$, $p_t$, $u_t$ & $\phi_t$ and sampled & historical data & historical data\\
		\hline
		variable: & $\psi_t$ & $u_t$ & - & -\\
		\hline
		depends on: & historical data & $(\phi_{t,t-1}, \psi_{t,t-1})$ & - & - \\
		\hline
	\end{tabular}
	\caption{Description of variables in the Wind Field Model and their dependence. \label{t:wfm}}
\end{table}

Throughout this article, we will define a \textit{hurricane event} as a tropical cyclone taking any form (e.g. tropical depression, tropical storm, hurricane) and denote the total lifetime of a hurricane event as a length of indexed time $T$ representing the total number of 6-hour time intervals passed since formation. For any given hurricane event, if we are given the \textit{track}, $(\phi_t,\psi_t)$, and central pressure timeseries, $p_t$, we may use (\ref{wfm}) to reconstruct the maximum windspeed $V(R_{\max},\phi_t,p_t,u_t)$ for all $t = 1,\dots, T$ where $t$ is the index number of 6-hour time intervals passed since formation. From (\ref{wfm}) we can see that extreme highs of the maximum windspeed occur for extreme lows of central pressure. Hence, it is important to accurately model the central pressure minima of a hurricane event in order to estimate longer year returns and rare threshold exceedances of maximum windspeeds. Furthermore, central pressure minima often occur at or near landfall so they are particularly important for estimating coastal risk.

We can use tools from extreme value theory to model extremes of a timeseries (e.g. minima or maxima). One well-known strategy is to approximate the set of maxima (or negative minima) taken over blocks of a fixed length $m$ of a set of independent and identically distributed random variables by the generalized extreme value distribution (GEV) given by,

\begin{equation}\label{gev}
G(x) = \exp\bigg[-\bigg\{1-\frac{k(x-\mu)}{\sigma}\bigg\}^{-1/k}\bigg]
\end{equation}
for $x~:~1+k(\frac{x-\mu}{\sigma})\ge0$ where $\mu$ is the location parameter, $\sigma$ is the scale parameter, and $k$ is the shape parameter that defines the tail behavior of $G$. Under certain regularity conditions, we may use maximum likelihood estimation of the parameters $\mu$, $\sigma$, and $k$ to fit the GEV to the block maxima (or negative minima) where  each parameter estimate is asymptotically normal provided $k > -0.5$ \cite{C}.

By a standard max-stable argument it is not necessary that the block length $m$ be fixed, as long as it is \textit{long enough} so that the maxima (or negative minima) can be modeled by its asymptotic GEV. By the same argument, the GEV that is fit to blocks of varying length is related to $G$ from (\ref{gev}), with different $\mu = \mu^*$ and $\sigma = \sigma^*$ parameters. A result of this max-stability property is that we may model the central pressure minima of hurricane events coming from historical records with varying lifetimes $T$. That is, we may model the negative central pressure minima
\begin{equation}
p_{\min} = \min_{t} p_t~~\text{for}~t = 1,\dots, T 
\end{equation}
by the GEV provided the lifetime of each hurricane event is long enough. We will refer to $t_{p_{\min}}$ as the time in the total lifetime of the hurricane event that the central pressure minimum is reached.

We can relax the requirement of strict independence for the GEV in (\ref{gev}) provided the time-series is weakly dependent and stationary, see for example~\cite[Chapter 3]{Leadbetter}
or \cite{exbook}. Using historical recordings from the HURDAT2 database, we find that the central pressure minima have the same dependence as in \cite{CC} on the lifetime $T$ and the latitude $\phi_{t_{p_{\min}}}$ where the central pressure minima occurs. Figure \ref{scatter} in the Appendix depicts scatter plots of the central pressure minima against the lifetime $T$ and latitude $\phi_{t_{p_{\min}}}$ for landfalling and non-landfalling hurricanes. 

We limit our model to hurricanes with lifetimes $T\ge25$ to ensure convergence of the negative central pressure minima to a GEV distribution. There are 642 hurricane events over the years 1851-2019 in the HURDAT2 database that satisfy this requirement. We use maximum likelihood estimation on the parameters of the stationary GEV model proposed in \cite{CC}. Although the central pressure minima $p_{\min}$ are sampled from independent hurricane events, they have some underlying dependence on both the lifetime $T$ and latitude $\phi_{t_{p_{\min}}}$ of the hurricane event which is accounted for in the location $\mu$ and scale $\sigma$ parameters of this stationary model. Figure \ref{SFIT} illustrates still poor fits for quantile plots of this stationary model. 

A natural question is whether there exists some time-dependence in the distribution of central pressure minima.

\subsection{A Nonstationary Model for Central Pressure Minima}

We investigate the time-dependence in the location $\mu$ and scale $\sigma$ parameters of central pressure minima for landfalling and nonlandfalling hurricanes. We perform an $F$-test for equal variance that indicates the variance of the central pressure minima for landfalling hurricanes has significantly changed ($p = 0.0037<0.05$) in the last 40 years. We obtain a similar result using a $T$-test for equal means of the central pressure minima for nonlandfalling hurricanes ($p = 0.0121<0.05$). Preliminary investigations into the time-dependence of the shape $k$ parameter showed no obvious trend, so it is taken as constant.

Motivated by the observed difference in the statistical parameters of the central pressure minima in the last 40 years, we now investigate the possibility of a time-dependent trend in the location and scale parameters of the stationary model proposed by \cite{CC}. This stationary model asserts a dependence of the location $\mu$ and scale $\sigma$ parameters on the lifetime $T$ of the hurricane and latitude $\phi_{t_{p_{\min}}}$ of the central pressure minima. We use this model as a basis for checking the time-dependence of the $\mu$ and $\sigma$ parameters in the GEV described by (\ref{gev}). We begin by performing maximum likelihood estimation of the all the coefficient parameters used in the stationary model where this estimation is performed on subsets of the 642 historical values of central pressure minima taken over moving time windows of 40 years with a timestep of 1 year. Our final result is a set of time-series representing the maximum likelihood values of the coefficient parameters in the stationary model. We then reconstruct the time-series of the location, $\mu(t_{\text{yr}})$, and scale, $\sigma(t_{\text{yr}})$, using the relationships described in the stationary model ($t_{\text{yr}} = \text{yr}-1851$) and the historical values of $T$ and $\phi_{t_{p_{\min}}}$. From now on, we will refer to the time-series $\mu(t_{\text{yr}})$ and $\sigma(t_{\text{yr}})$ as the \textit{location time-series} and \textit{scale time-series}, respectively, to differentiate between the other time-series in this investigation.

Unreliable maximum likelihood estimates of the location and scale parameters in the years 1851-1960 are due to low numbers of recorded hurricane events. Nevertheless, continuous time-dependent trends are noticeable after 1960 for parameters in both the landfalling and nonlandfalling case.

A Mann-Kendall test for trend is performed on each of the location and scale time-series constructed as described above from the 642 historical recordings of central pressure minima. We remark that there are 300 location time-series in the landfalling case and 342 location and scale time-series in the nonlandfalling case. This is because our estimated location (similarly, location and scale) parameter(s) depend on the lifetime $T$ and the location $\phi_{t_{p_{\min}}}$ of minimum central pressure of the hurricane where we have 300 (similarly, 342) historical recordings of such lifetimes and locations. This is in contrast to the scale parameter of landfalling hurricanes which does not have a dependence on $T$ or $\phi_{t_{p_{\min}}}$ and, as a consequence, results in a single scale time-series.

\begin{figure}
	\centering
	\includegraphics[width=0.5\textwidth]{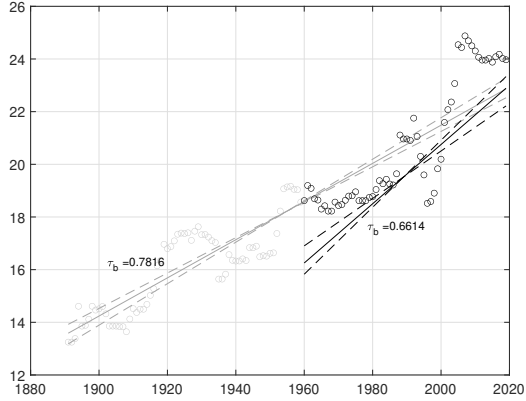}
	\caption{Timeseries of parameter $\sigma(t_{\text{yr}}) = \sigma_0(t_{\text{yr}})$ coming from the stationary model for $-p_{\min}$ in hPa of landfalling hurricanes constructed from likelihood estimates. The value $\tau_b$ is the Kendall correlation coefficient. \label{lf_ts}}
\end{figure}

\begin{figure}
\begin{minipage}{0.5\textwidth}
	\centering
		\includegraphics[width=\textwidth]{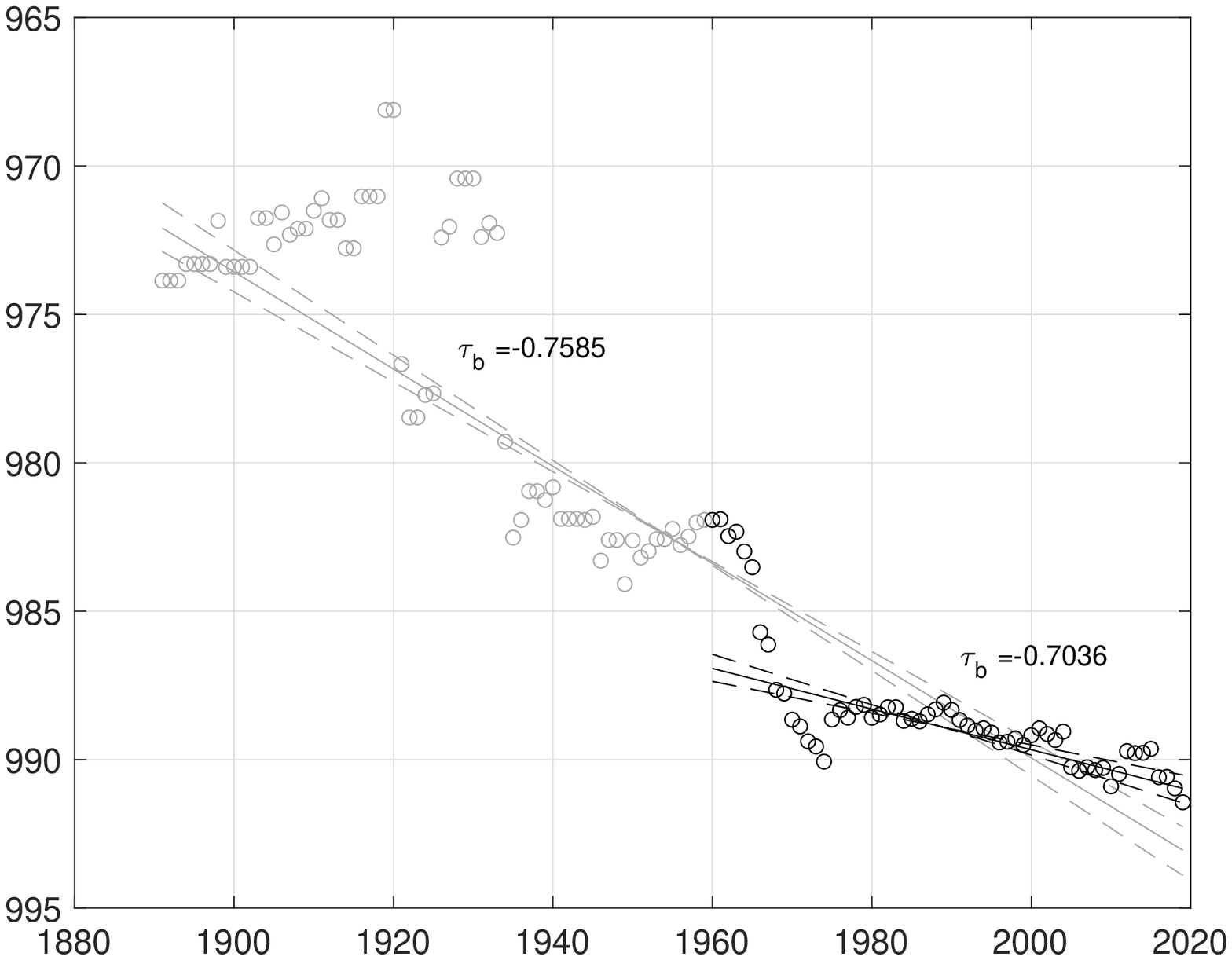}
		\caption*{(a)}
\end{minipage}
\begin{minipage}{0.5\textwidth}
	\centering
		\includegraphics[width=\textwidth]{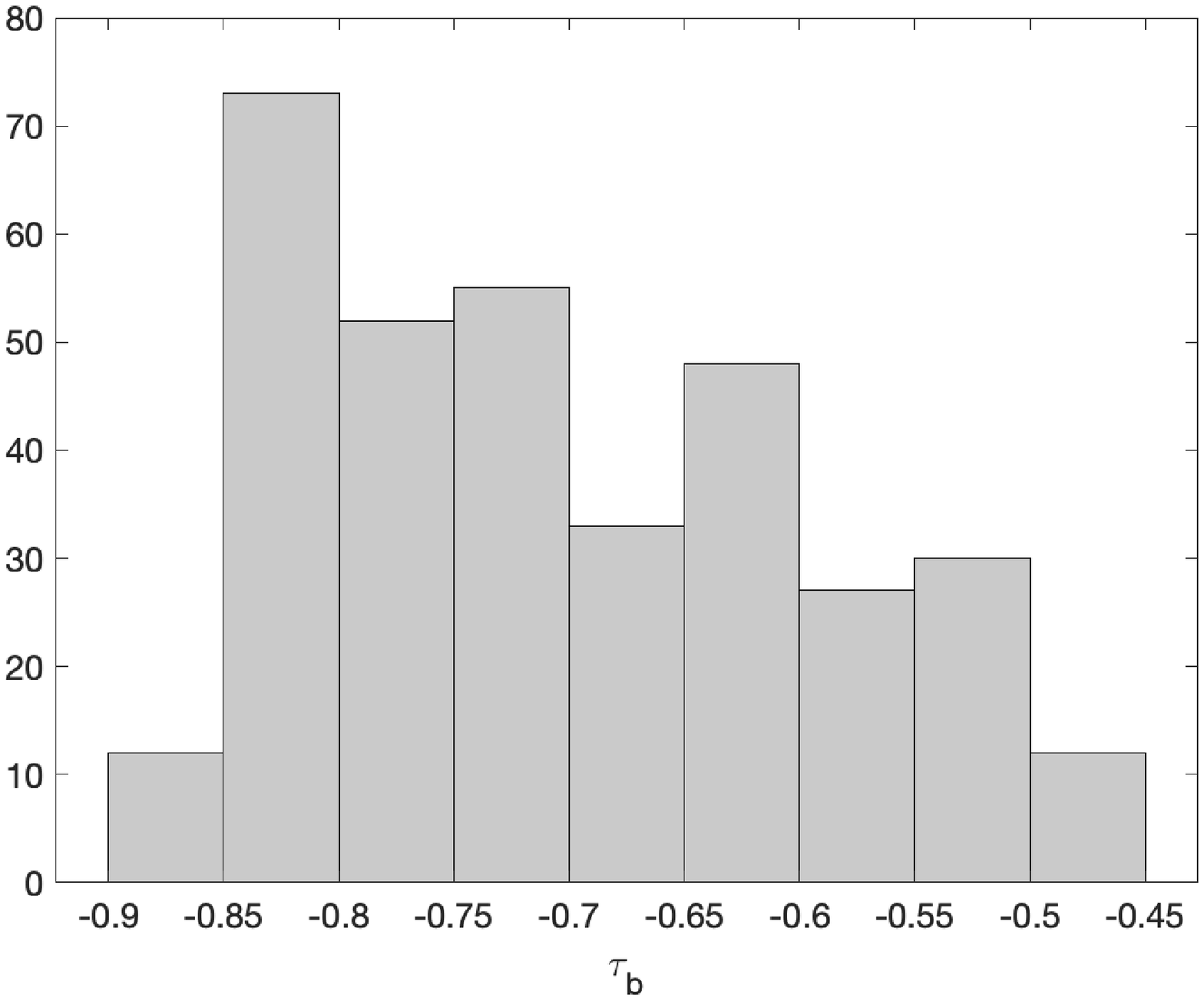}
		\caption*{(b)}
\end{minipage}
	\caption{(a) An example time-series of parameter $\mu(t_{\text{yr}}) = \mu_0(t_{\text{yr}})+\mu_1(t_{\text{yr}})*\log(T)$, for a single sampled $T$. The value $\tau_b$ is the Kendall correlation coefficient. The $y$-axis represents the end year of the 20-year time window chosen for likelihood parameter estimation (b) $\tau_b$ for all time-series of $\mu(t_{\text{yr}})$ illustrating that all time-series have a statistically significant negative Kendall correlation coefficient. \label{nlf_ts}}
\end{figure}

To determine whether a trend is reliable, we perform the Mann-Kendall test for trend on all the location and scale time-series. We find a positive statistically significant trend for the scale parameter in the landfalling case and a negative statistically significant trend for all time-series of the location parameter in the nonlandfalling case. We do not find clear evidence for a reliable trend in the location time-series for the landfalling case or the scale time-series for the nonlandfalling case. Although, some of the 300 (similarly, 342) location (similarly, scale) time-series do exhibit a significant trend. All trends are determined at the $\alpha = 0.05$ significance level. The Kendall correlation coefficient is estimated for all years and for years from 1960-2020, for comparison. For an illustration of the trend results and parameter time-series see Figures \ref{lf_ts} and \ref{nlf_ts}. 


From these results, we propose the following nonstationary model ($t_{\text{yr}} = \text{yr}-1851$ is the yearly index),
 \begin{equation}\label{ns_lf}
 \mu = \mu_0+\mu_1\log(T)+\mu_2\phi_{t_{\text{p}_{\min}}},~~~~\sigma = \sigma_0+\sigma_1t_{\text{yr}},~~~~k = k_0
 \end{equation}
 for landfalling hurricanes, and
 \begin{equation}\label{ns_nlf}
 \mu = \mu_0+\mu_1\log(T)+\mu_2\log(t_{\text{yr}}),~~~~\sigma = \sigma_0+\sigma_1\phi_{t_{\text{p}_{\min}}},~~~~k = k_0
 \end{equation} 
for nonlandfalling hurricanes. Maximum likelihood estimates and standard errors estimated from the information matrix are provided in Table \ref{tdgev}. We find using the likelihood ratio test that our revised nonstationary model for central pressure minima offers a statistically significant better fit to the data than the stationary model with test statistics well beyond $\mathcal{L}_{0.05,1} = 3.84$, the statistic corresponding to the $\alpha=0.05$ significance level with 1 degree of freedom (see $\mathcal{L}$ in Table \ref{tdgev}). 

\begin{table}[ht]
	\centering
	\begin{tabular}{|lllll|}
		\hline
		Type & $\mu_0$(se) & $\mu_1$(se) & $\mu_2$(se) & \\
		\hline
		Landfalling  & -1078.97(14.48) & 32.87(3.60) & -0.52(0.16) & \\
		Nonlandfalling & -1027.13(13.65) & 27.40(2.76) & -11.47(1.84)*& \\
		\hline
		Type & $\sigma_0$(se) & $\sigma_1$(se) & $k$(se) & $\mathcal{L}$\\
		\hline
		Landfalling & 12.47(1.86) & 0.07(0.02)* & -0.13(0.04) & 15.62\\
		Nonlandfalling & 20.48(2.04) & -0.16(0.06) & -0.13(0.04) & 30.77\\
		\hline		
	\end{tabular}
	\caption{Maximum likelihood estimates of the parameters in the nonstationary generalized extreme value distribution model for $-p_{\min}$. Time-dependent parameters are marked with *. Likelihood ratio test statistics for our revised nonstationary model of $-p_{\min}$ against the stationary model are indicated by $\mathcal{L}$. \label{tdgev}}
\end{table}

\subsection{Poisson Returns of Hurricane Events and a Nonstationary Rate Parameter}

We investigate a Poisson model for yearly returns of hurricane events where the number of expected yearly hurricane events is increasing over time. 

Returns of extreme hurricane events, such as low central pressure minima or high maximum windspeeds, are often reported in terms of an $n$-year return. In order to interpret returns in this way, our model must consider how often a hurricane event occurs in a given year. Classically, it is expected that a rare event, such as a hurricane, is modeled by a Poisson distribution given by,
\begin{equation}\label{poisson}
P(X=\mathscr{K}) = \frac{\lambda^{\mathscr{K}}e^{-\lambda t}}{\mathscr{K}!}
\end{equation}
where $\lambda = r/t$ is the rate parameter estimated as the number of events $r$ in a given time $t$. 

\begin{figure}
	\centering
	\includegraphics[scale=0.12]{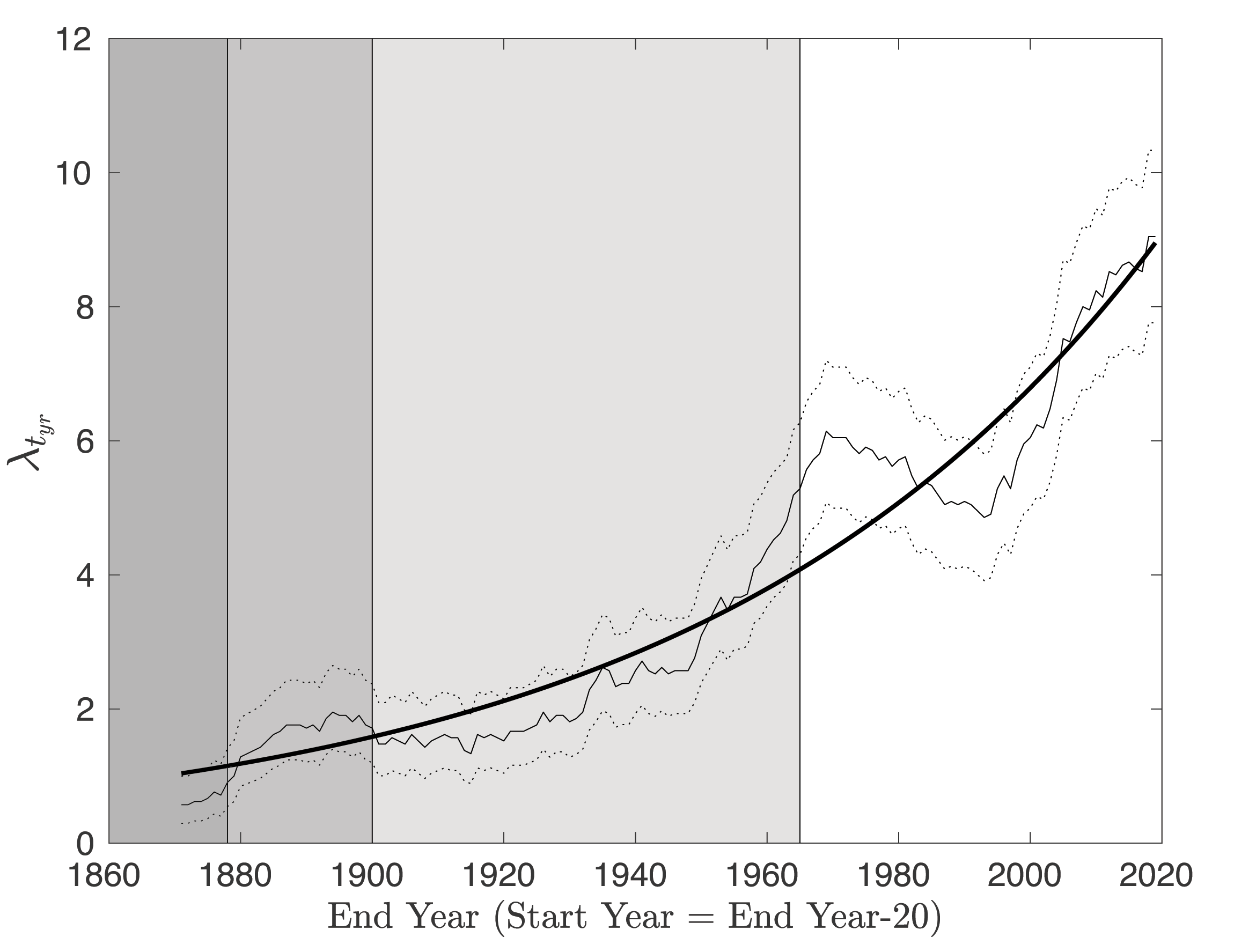}
	\caption{Likelihood estimate of the Poisson parameter for yearly hurricane event rates with lifetimes greater than 6.25 days. Estimates are taken over 20 year moving time windows. Standard errors are marked with dotted lines. Fitted exponential model is represented by a thick line. Grayed areas correspond to those in \cite{V3}: (1) 1878 - year when the U.S. Signal Corps began cataloging all Atlantic hurricanes (2) 1900 - year when the U.S. Coast was sufficiently well-populated for monitoring (3) modern-era with appropriate ship density. \label{POISSON_PARAM}}
\end{figure}

Under the assumption of stationarity, the authors in \cite{CC} estimate a fixed rate parameter, $\lambda=5.45$ hurricane events per year, as the average number of returns of a hurricane in a given year over the years 1965-1994. With more data available in the HURDAT2 database, we are able to estimate the time-dependent yearly rate parameter $\lambda_{t_{\text{yr}}}$ over 20 year sliding windows from 1851-2019. We refer the reader to Figure \ref{POISSON_PARAM} for an illustration of the estimated yearly rate parameter. 
\begin{remark}
	We use only hurricane events that have lifetimes $T\ge 25$ in this analysis to estimate the time-dependent rate parameter $\lambda_{t_{\text{yr}}}$; however, all hurricane events are observed at a further increased rate of $\lambda_{\text{2019}}\approx 16$ which corresponds to current NOAA estimates.
\end{remark}

We use likelihood estimation to fit an exponential to the time-dependent Poisson rate parameter $\lambda_{t_{\text{yr}}}$. In particular, our model is given as,
\begin{equation}\label{rt}
\lambda(t_{\text{yr}}) = ae^{bt_{\text{yr}}}.
\end{equation}
Likelihood estimates and confidence intervals of $a$ and $b$ can be found in Table \ref{rate}. Our model for hurricane returns does not differentiate between landfalling and nonlandfalling hurricane events due to the nature of the simulation in the final section. This is because tracks of a simulated hurricane are generated by randomly sampling a historical track and adding noise. To compare our results against current literature, we separate the discussion of returns of hurricane events for landfalling and nonlandfalling hurricanes in the paragraphs below.

There is some debate on whether the average number of hurricane events is increasing generally; some literature suggests that low ship density is the underlying cause for the low number of recorded hurricanes for years up to 1965 \cite{L,V2}, while others report significant increases in frequency after the late 1980s \cite{V}. When averaging yearly frequency over moving time windows, the authors in \cite{V} report a small nominally positive upward trend post 1878. The work of \cite{L} finds an increase in the occurrence of short lifetime hurricanes only, leading the authors to conclude ship density as a plausible cause for the observed trend. It is important to note that the literature described here uses the retired HURDAT database for their analyses rather than the HURDAT2 database used in this investigation; however, this certainly does not rule out the possibility of historically unrecorded storms in the updated database. There is active research on the frequency of hurricane events recorded in the HURDAT2 database where an observed late-20th century trend is attributed to a possible unusually low minima in the 1980s. Currently, this research is only available in preprint form in \cite{V3}.

We find an increasing trend in frequency of hurricane events longer than 6.25 days using the Poisson rate parameter, which differs from the results in \cite{L}. This trend holds even into the modern era (post 1965) where ship density is expected to remain steady. One explanation for this difference could be our use of a Poisson rate estimate over a moving average. Rate estimates expect that an increase in the mean results in an increase in the variance. This phenomenon is observed in the raw data. In the case of a moving average estimate this increase in variance can cause statistical tests of the mean difference to be near zero due to large standard errors. We also do not separate hurricane events by windspeed where differences in trend have been reported \cite{V}. Since we limit our investigation to hurricanes with lifespans longer than 6.25 days, our findings may also be a result of some underlying increase in the lifespan of hurricane events as a whole.  Finally, using the yearly estimates from our model for the rate over 1965-94 we find that the average is identical to past literature \cite{CC} which provides some reasonable benchmark.

An argument could be made that this increase in the total number of observed hurricane events post-1965 comes from our ability to more readily observe nonlandfalling hurricanes. However, an increase in the Poisson rate parameter is also observed for strictly landfalling hurricane events of lifetimes longer than 6.25 days;  however, this rate parameter follows a similar pattern (with a low minima in the 1980s) to that of \cite{V3} with a slight increase in the current peak compared to that of 1965. We refer to the Appendix figure \ref{POISSON_PARAM_LF} for an illustration of the estimated Poisson rate parameter for landfalling hurricane events and figure \ref{HURR_NUM} for the raw data of total yearly observed hurricane events, both with lifetimes longer than 6.25 days.

\begin{table}[ht]
	\centering
	\begin{tabular}{|ll|}
		\hline
		$a$ (ci) & $b$ (ci)\\
		\hline
		1.024 (0.925, 1.123) & 0.015 (0.014, 0.015)\\
		\hline
	\end{tabular}
	\caption{Maximum likelihood estimates of the exponential model for the time-dependent Poisson parameter $\lambda_{t_{\text{yr}}}$.\label{rate}}
\end{table}

\subsection{Verification of the Nonstationary Model for Central Pressure Minima}

We use a combination of established statistical methods to illustrate the reliability of our nonstationary model at predicting returns of central pressure minima. 

To test the reliability of our model to accurately predict the distribution of central pressure minima without updating, we break the HURDAT2 database up into a \textit{training} set which we will use to simulate hurricanes from the model and \textit{test} set which we will use to compare risk probability outcomes estimated from the training set against the 'true' probabilities. Our training set will be defined as the set of all years in our dataset minus the number of years $n$ used to obtain the $n$-year returns and our test set will be the $n$ last years in our dataset. For example, if we are interested in finding the $50$-year returns, our training set would be defined as the set of all hurricanes occurring between 1851-1970 and our test set would be the set of all hurricanes occurring between 1971-2019.

Under the assumption that our negative central pressure minima follow some generalized extreme value distribution, \cite[Section 6.2.3]{C} suggests the use of a sequence of standardized variables $z_{t_{\text{yr}}}$ defined for our purposes by,
\begin{equation}\label{stand_z}
z_{t_{\text{yr}}} = \frac{1}{k} \log\bigg\{1+k\bigg(\frac{-p_{\min}(t_{\text{yr}})-\mu(t_{\text{yr}})}{\sigma(t_{\text{yr}})}\bigg)\bigg\}
\end{equation}
each having a standard Gumbel distribution,
\begin{equation}\label{gumbel}
P(z_{t_{\text{yr}}}\le z) = \exp\{-e^{-z}\},~~~~z\in\mathbb{R}.
\end{equation}
The advantage of using this sequence is that the 'true' quantile plots of the observed and standardized $-p_{\min}(t_{\text{yr}})$ in the test set can be made with reference to the distribution for the simulated and standardized $-p_{\min}(t_{\text{yr}})$ from the training set. 

\begin{figure}
	\centering
	\begin{minipage}{0.45\textwidth}
		\centering
		\includegraphics[scale=0.35]{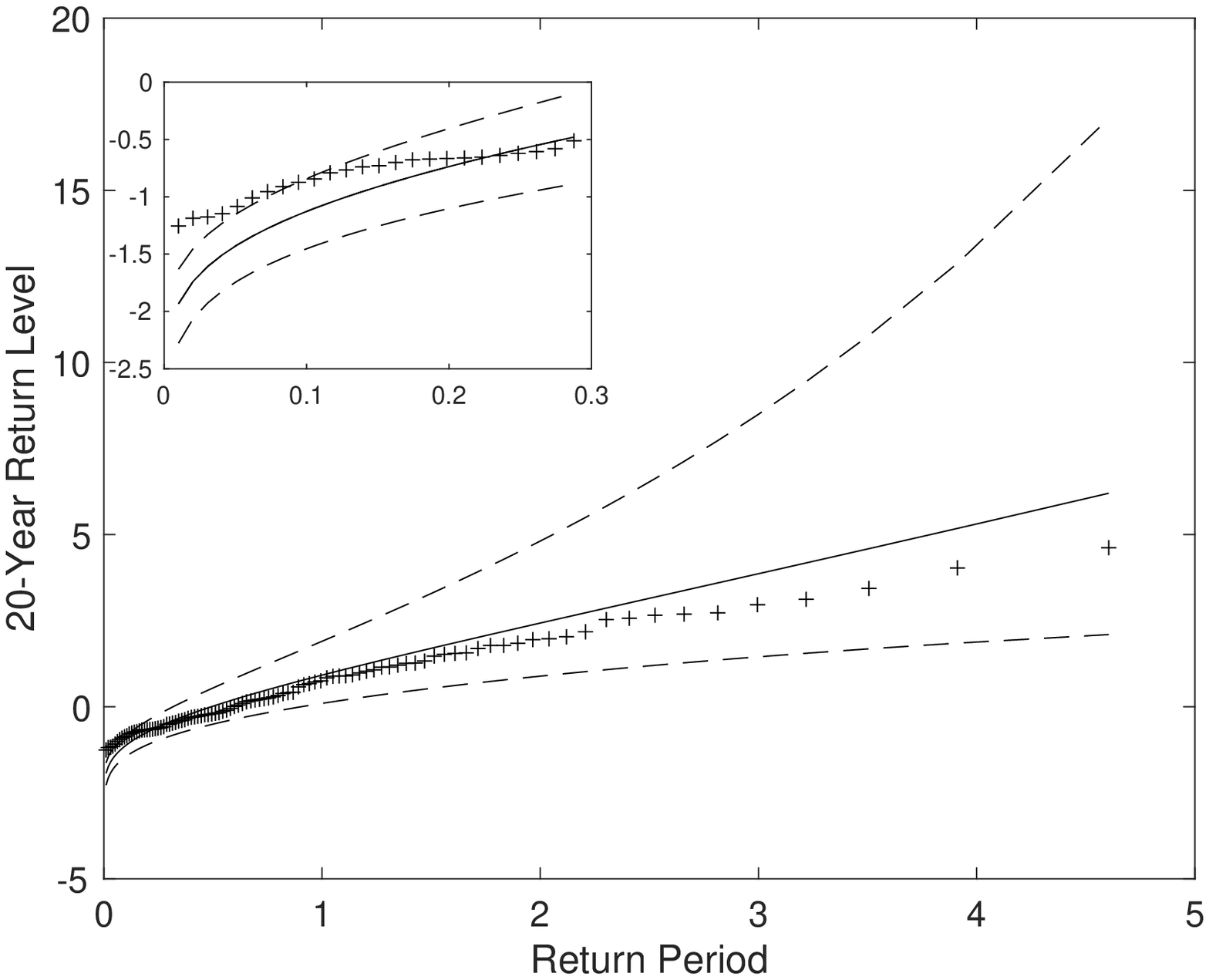}
		\caption*{(a)}
	\end{minipage}
	\begin{minipage}{0.45\textwidth}
		\centering
		\includegraphics[scale=0.35]{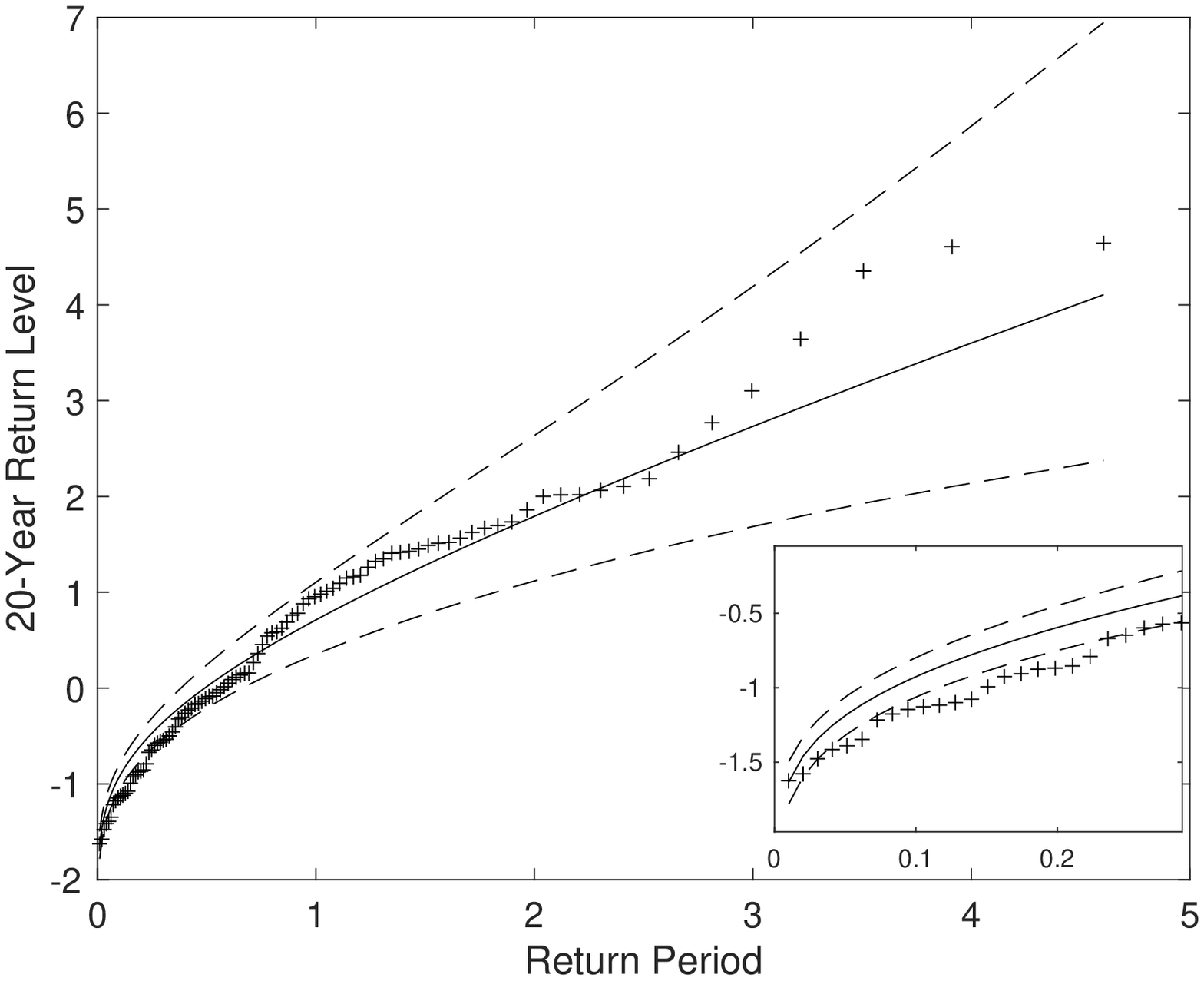}
		\caption*{(b)}
	\end{minipage}
	\caption{$20$-year return-levels of the standardized $-p_{\min}$ coming from $(\ref{stand_z})$ for (a) nonlandfalling and (b) landfalling hurricanes. Solid lines and dashed lines represent the model and 95\% confidence intervals approximated from the training set over the years 1851-2000. Symbol (+) indicates the true return-levels calculated from the test set over the years 2001-2019. Return periods and return levels here are based on $(\ref{stand_z})$ are non-dimensional and expected to follow the Gumbel distribution $(\ref{gumbel})$. \label{qq_20}}
\end{figure}

\begin{figure}
	\centering
	\begin{minipage}{0.45\textwidth}
		\centering
		\includegraphics[scale=0.35]{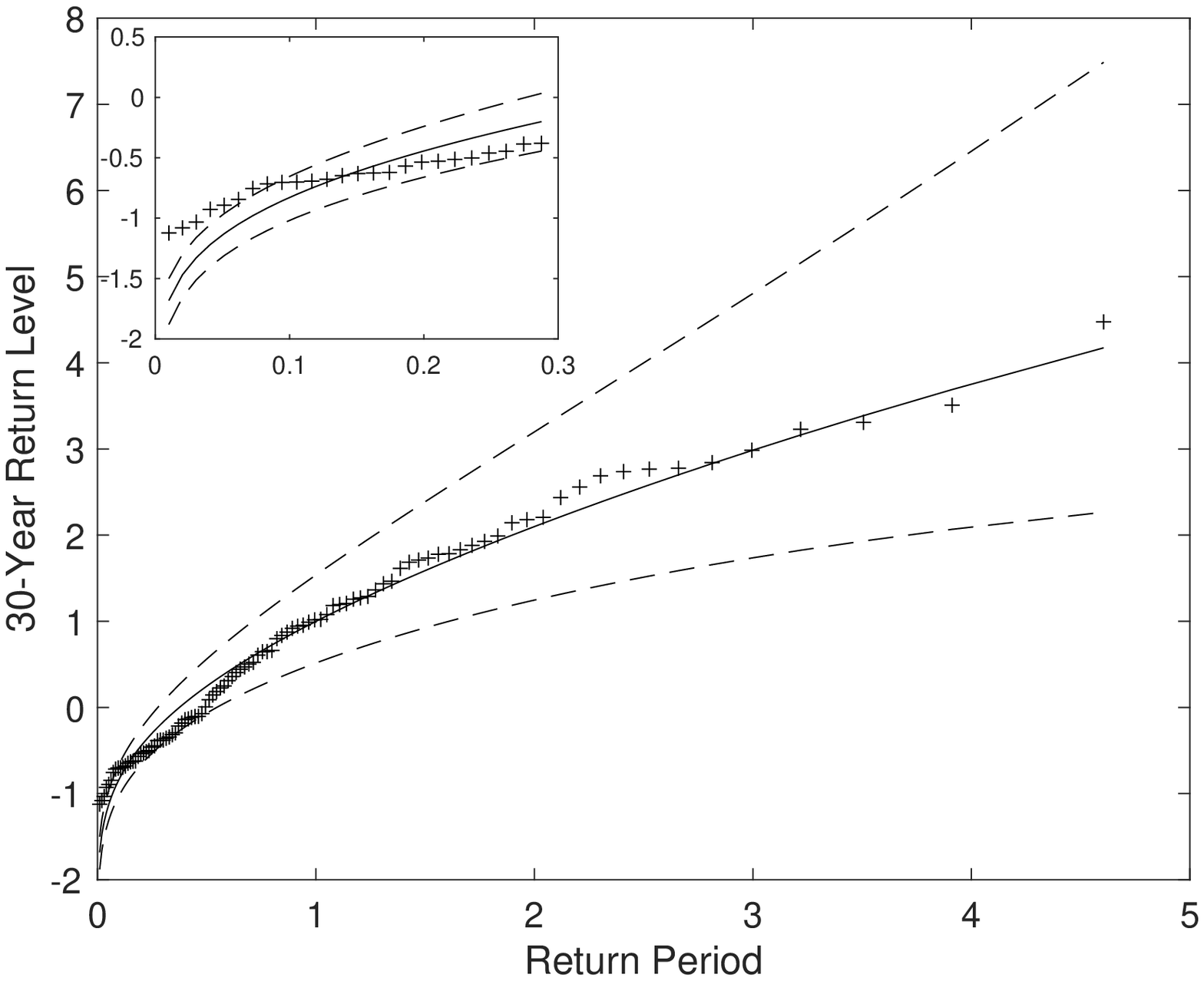}
		\caption*{(a)}
	\end{minipage}
	\begin{minipage}{0.45\textwidth}
		\centering
		\includegraphics[scale=0.35]{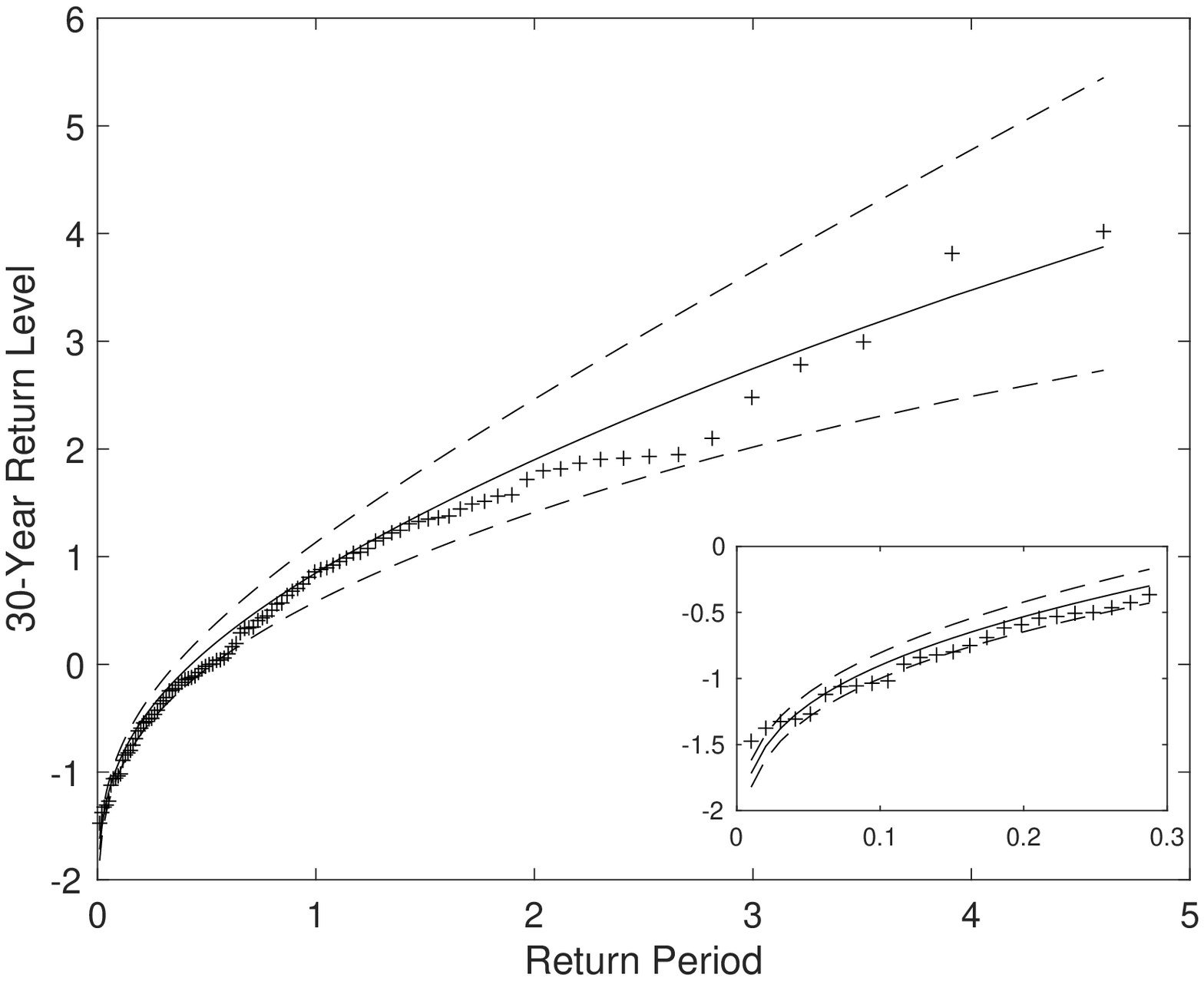}
		\caption*{(b)}
	\end{minipage}
	\caption{$30$-year return-levels of the standardized $-p_{\min}$ coming from $(\ref{stand_z})$ for (a) nonlandfalling and (b) landfalling hurricanes. Solid lines and dashed lines represent the model and 95\% confidence intervals approximated from the training set over the years 1851-1990. Symbol (+) indicates the true return-levels calculated from the test set over the years 1991-2019. Return periods and return levels here are based on $(\ref{stand_z})$ are non-dimensional and expected to follow the Gumbel distribution $(\ref{gumbel})$. \label{qq_30}}
\end{figure}

\begin{figure}
	\centering
	\begin{minipage}{0.45\textwidth}
		\centering
		\includegraphics[scale=0.35]{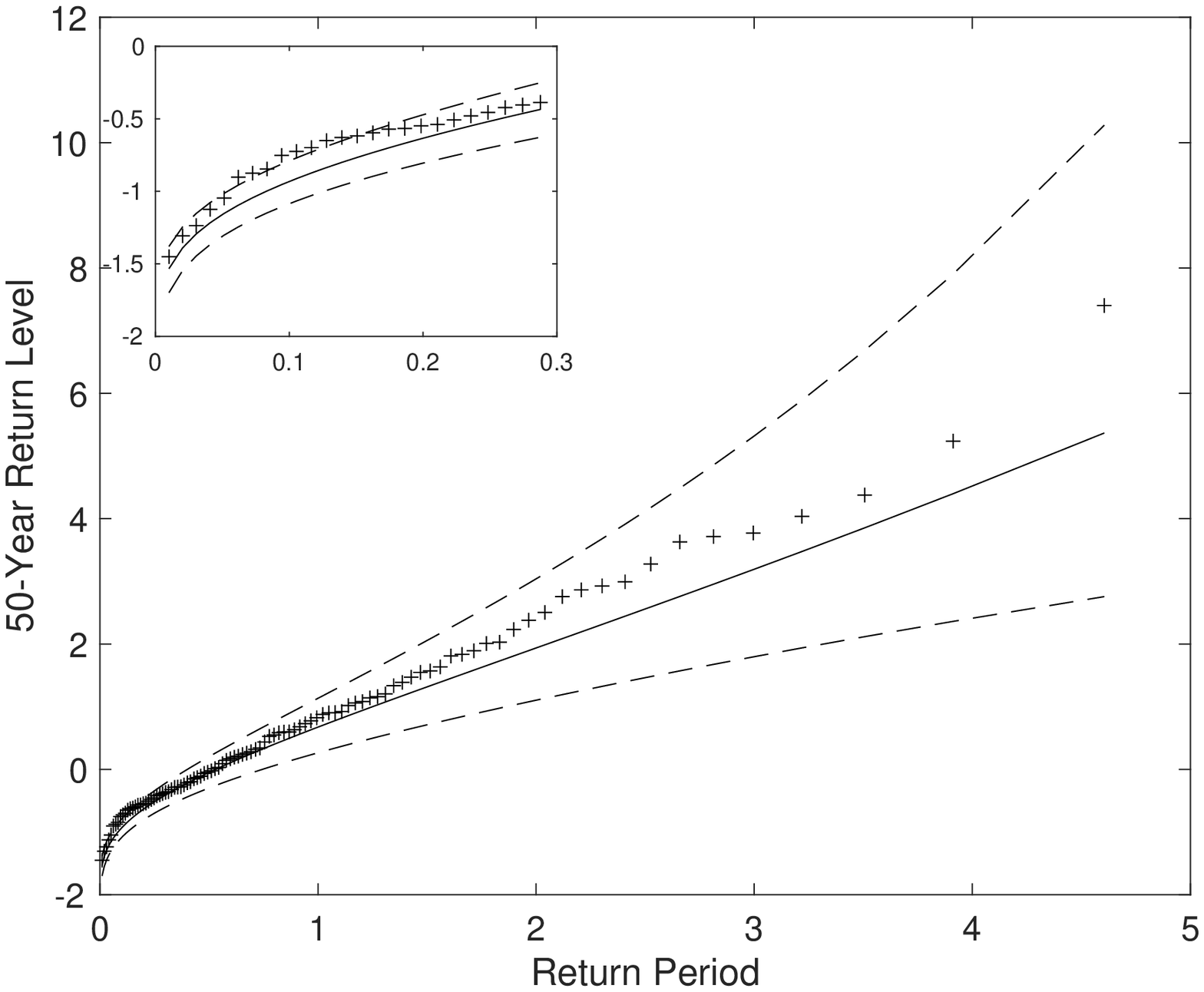}
		\caption*{(a)}
	\end{minipage}
	\begin{minipage}{0.45\textwidth}
		\centering
		\includegraphics[scale=0.35]{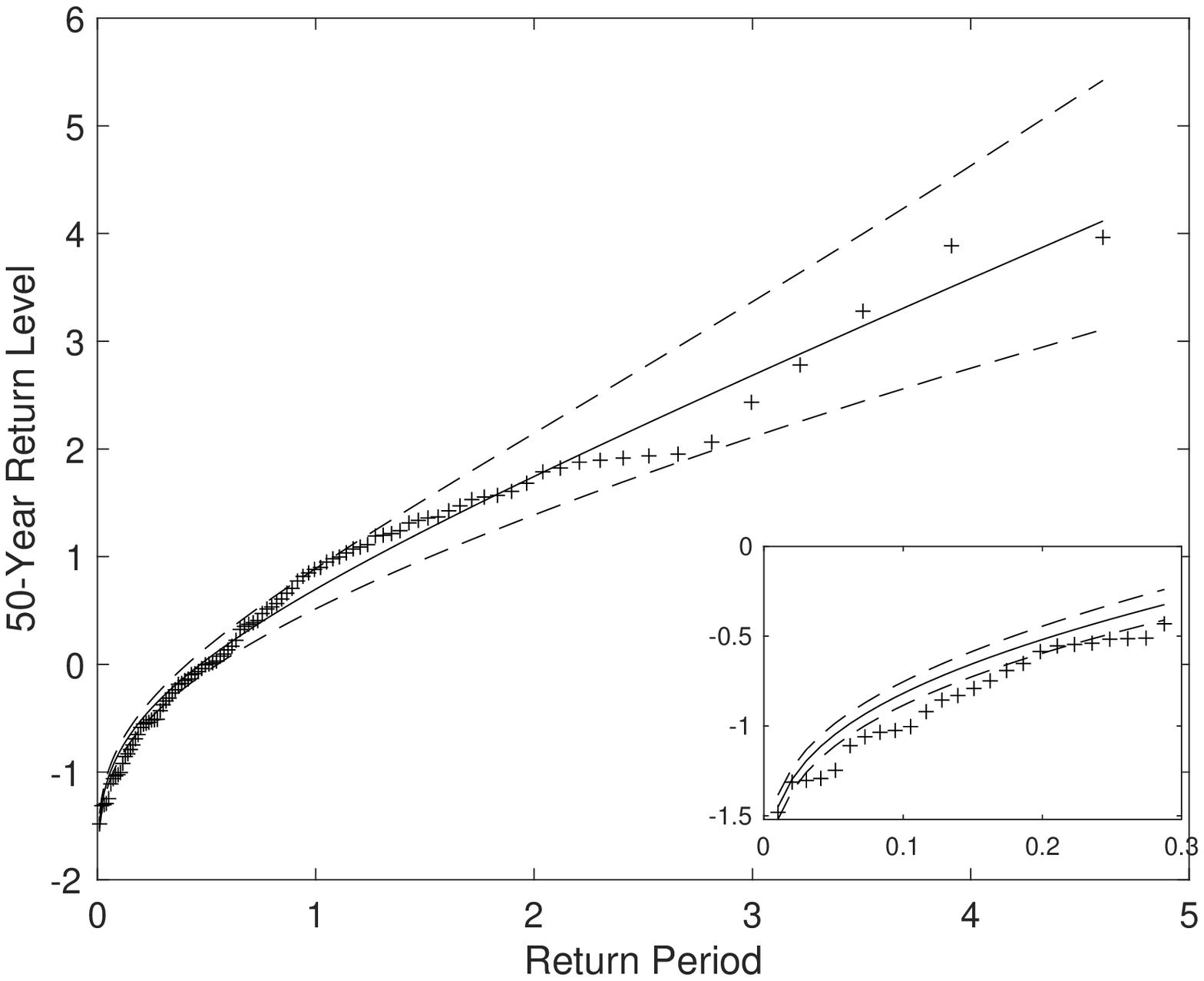}
		\caption*{(b)}
	\end{minipage}
	\caption{$50$-year return-levels of the standardized $-p_{\min}$ coming from $(\ref{stand_z})$ for (a) nonlandfalling and (b) landfalling hurricanes. Solid lines and dashed lines represent the model and 95\% confidence intervals approximated from the training set over the years 1851-1970. Symbol (+) indicates the true return-levels calculated from the test set over the years 1971-2019. Return periods and return levels here are based on $(\ref{stand_z})$ are non-dimensional and expected to follow the Gumbel distribution $(\ref{gumbel})$. \label{qq_50}}
\end{figure}

We generate data to model negative central pressure minima $n$-year returns for the years in the test set using the 1.) parameter likelihoods of $\mu(t_{\text{yr}})$, $\sigma(t_{\text{yr}})$ and $k$ defined by the model in (\ref{ns_lf}) and (\ref{ns_nlf}) estimated from the training set and 2.) the appropriate rate parameters defined by (\ref{rt}) to compute returns of hurricane events using (\ref{poisson}) where $t_{\text{yr}}$ indices are chosen to correspond to those of the test set. We use this data to compute the model standardized quantile plots for $20$-, $30$-, and $50$-year return periods for both landfalling and nonlandfalling hurricanes. Figures \ref{qq_20}, \ref{qq_30}, and \ref{qq_50} show model results against the actual data in the test set. Not surprisingly, better approximations for both the landfalling and nonlandfalling case are made for shorter $n$-year returns; however, estimates for $50$-year returns still fall reasonably within the 95\% confidence interval of the model estimated from the information matrix.

Using the standardized negative central pressure minima allows us to estimate the accuracy of the nonstationary model against true data; however, it does not provide us with a complete way of interpreting the $n$-year returns. At best, we are able to fix a year index $t_{\text{yr}}$ and state the probability of the negative central pressure minima being above a certain threshold in that given year. Most risk analysis involves directly computing $n$-year return-levels where a new definition needs to be introduced in the nonstationary setting. We discuss this in detail in the next section.

\section{Application of Methodology for Coastal Windspeed Risk}

\subsection{Time-dependent Returns of High Maximum Windspeeds}
We discuss a definition for time-dependent $n$-year return-levels of maximum windspeeds. 

Return-level is often used in risk analysis to communicate the threshold that we are expected to exceed in a give amount of time. For example, we may ask what is the maximum value of the windspeed that we are expected to exceed in $n$ years. When accounting for nonstationary effects, such as those brought on by climate change, the probability of observing values above or below a threshold varies over time so that terms like return-level no longer make physical sense.


\cite[Section 4.2]{CO} introduces the idea of extending  the definition of the $n$-year return-level to the nonstationary case by taking the threshold where the expected number of exceedances in $n$ years is 1 to the non-stationary case. In the context of nonstationary windspeed prediction this would be equivalent to solving for $r_n$ in (\ref{returns}),
\begin{equation}\label{returns}
1 = \sum_{t_{\text{yr}} = 1}^n (1-F_{t_{\text{yr}}}(r_n))
\end{equation}
where $r_n$ is the $n$-year return-level beginning with year $t_{\text{yr}} = 1$ and ending with year $t_{\text{yr}}=n$ and $F_{t_{\text{yr}}}$ is the unknown indexed yearly cumulative distribution function of maximum windspeed. For example, if we are interested in finding the $50$-year return-level $r_{50}$ of windspeed, (\ref{returns}) would become, 
\begin{equation}
1 = \sum_{t_{\text{yr}} = 1}^{50} (1-F_{t_{\text{yr}}}(r_{50})) = P_{1}(\text{ws}>r_{50})+\dots+P_{50}(\text{ws}>r_{50})
\end{equation}
The corresponding $n$-year return-level can be numerically estimated for future years by extrapolating the trend in the model and approximating $r_n$ by calculating the $1-\frac{1}{n}$ quantile of the equal weight mixed probability density function of windspeed occurring over $t_{\text{yr}} = 1,\dots,n$ years given by,
\begin{equation}\label{mixed}
f(x;t_{1},\dots,t_{n}) = \sum_{t_{\text{yr}} = 1}^{n} f_{t_{\text{yr}}}(x).
\end{equation}
where $f_{t_{\text{yr}}}$ is the unknown and numerically approximated probability density function of the windspeed corresponding to the yearly time index $t_{\text{yr}}$. In fact, the definition in (\ref{mixed}) has also been used to model regional returns of extremes where $f(x;\ell_1,\dots,\ell_n)$ varies by location $\ell_i$ instead of time \cite{CK}.

\subsection{A Simulation to Estimate Maximum Windspeed Risk Along the US North Atlantic Coast}
We run a simulation using the adaptations described in earlier sections to estimate high maximum windspeed risk for specified regions along the US North Atlantic coast.

From the Wind Field Model described in (\ref{wfm}), we observe that returns of low central pressure minima have a large and direct effect on returns of high maximum windspeeds. This relationship makes appropriately modeling central pressure minima vital when considering returns of extreme windspeeds along the coast. However, it is not enough to know the central pressure minima to estimate coastal windspeed risk. This is because maximum windspeeds for a coastal region depend, among other things, on the translational velocity of the hurricane, the location at which landfall occurs, and whether the central pressure minima is achieved at landfall. 

\begin{figure}
	\centering
	\begin{minipage}{0.48\textwidth}
		\includegraphics[width=\textwidth]{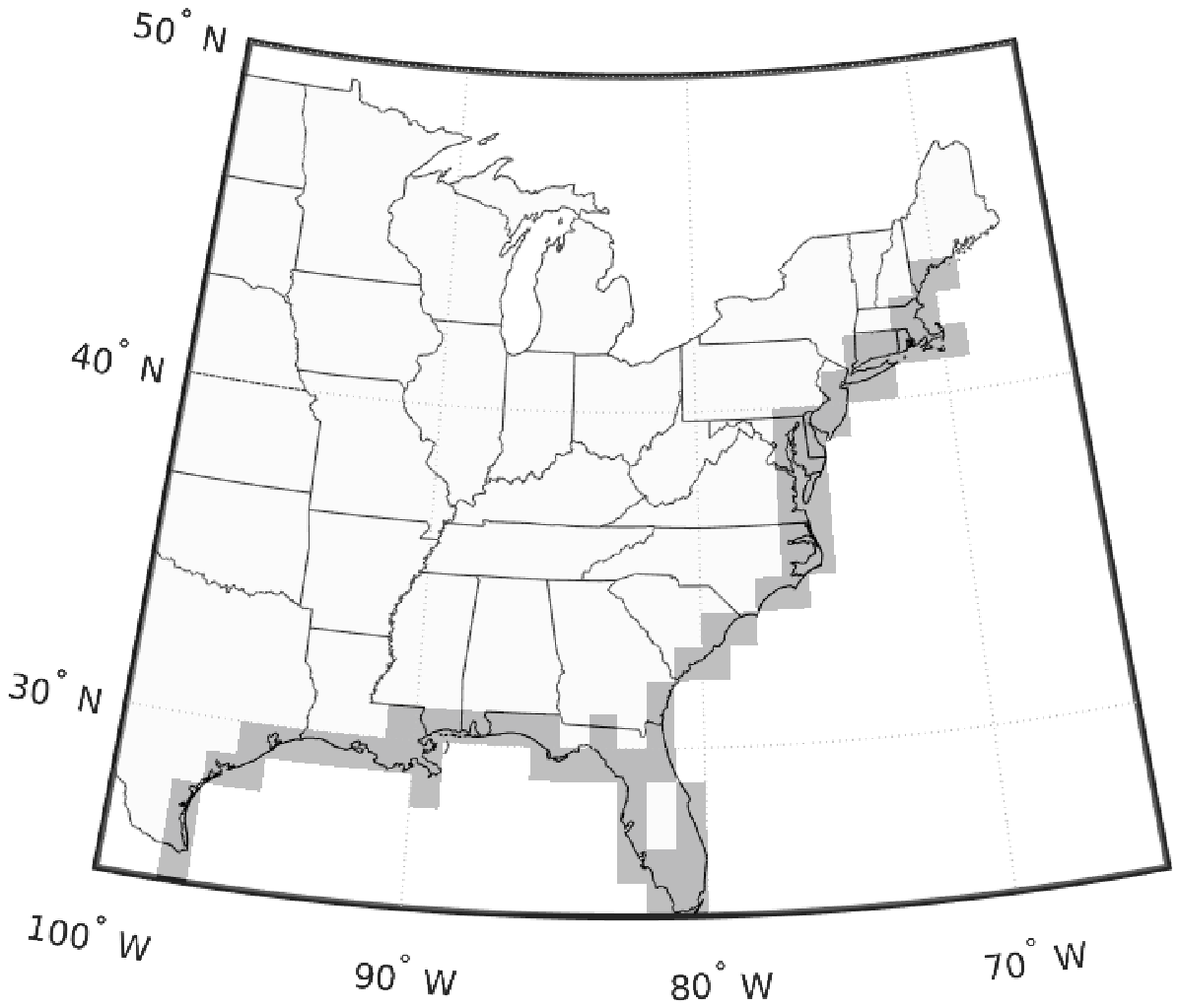}
		\caption*{(a)}
	\end{minipage}	
	\begin{minipage}{0.48\textwidth}
		\includegraphics[width=\textwidth]{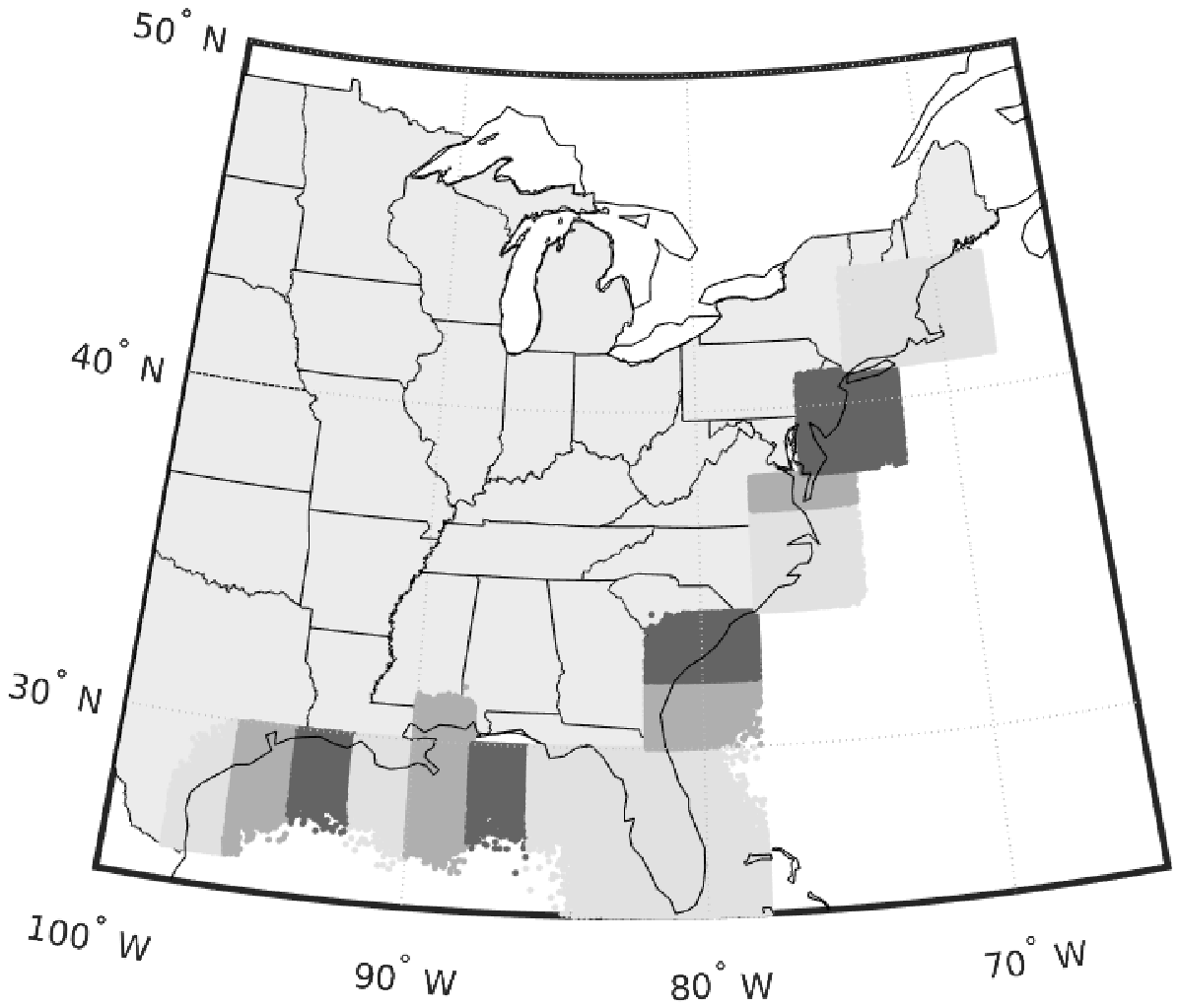}
		\caption*{(b)}
	\end{minipage}
	\caption{(a) Coarse grid representing the coastal line. (b) Simulated hurricane locations along the coastal line. Different regions are indicated in gray-scale.\label{coast}}
\end{figure}

We now consider a more complex hurricane simulation to estimate the unknown distribution described in (\ref{mixed}) of maximum windspeeds for a particular coastal location with the adaptations described in this investigation. The simulation is outlined in the Appendix \ref{simu}; however, we refer the reader to the original literature \cite{CC} for a detailed description. In essence, the process described in the Appendix \ref{simu} simulates a series of hurricane events for a given year by sampling the number of events to occur and the random variables used in the Wind Field Model represented by (\ref{wfm}) at each time $t$ along a simulated hurricane track. Once all hurricane events for a set of years have been simulated, we sample the windspeed for each simulated hurricane landing along a specified coastline to form the unknown distribution described in (\ref{mixed}). The North Atlantic coastline is first approximated by a coarse grid, illustrated in figure \ref{coast}, then divided into coastal regions: N-Texas, S-Texas, W-Louisiana, E-Louisiana, Mississippi, Alabama-Florida, Florida, Florida-Georgia, South Carolina, North Carolina, Virginia, Maryland-New Jersey, and Connecticut-Massachusetts-New Hampshire. A simulated hurricane is said to be "on the coast" if the eye of the hurricane is within 2 degrees of the coastal line.

To estimate the $n$-year return-levels for regions along the coast, we must numerically approximate the probability distribution function of windspeeds described in (\ref{mixed}). Then the $n$-year return-level is simply the $1-\frac{1}{n}$ quantile of the combined frequency distribution of maximum windspeed data for each coastal region. We do this for $20$-, $30$- and $50$-year return-levels for each region taken along the coast by generating 20, 30 and 50 years of data (that is, 2020-2040, 2020-2050, and 2020-2070) for $N = 1,000$
trials and estimating the $0.95$, $0.97$, and $0.98$ quantiles, respectively. 

It is reasonable to assume that each likelihood parameter in our simulation of maximum windspeeds (there are several), $\theta$, has reached its asymptotic normal distribution $\mathcal{N}(\hat{\theta},s_{\theta})$ with mean $\hat{\theta}$ equal to the maximum likelihood estimate of the parameter $\theta$ and standard deviation given by the standard error $s_{\theta}$ approximated from the Hessian. We can be confident that the true population distribution of maximum windspeeds, which the model is meant to represent, falls within some combination of these parameters; each coming from their corresponding distribution $\mathcal{N}(\hat{\theta},s_{\theta})$. 

\begin{figure}
	\centering
	\begin{minipage}{0.9\textwidth}
		\includegraphics[width=\textwidth]{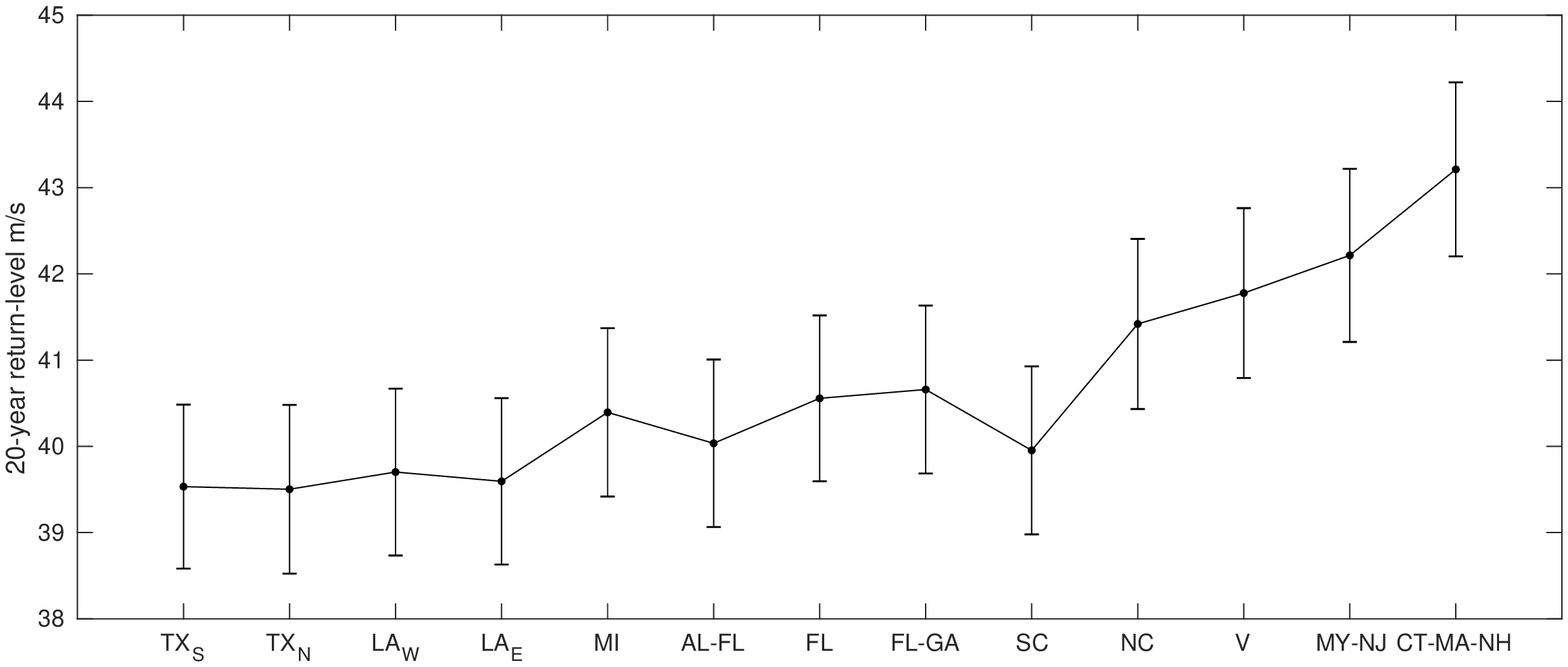}
	\end{minipage}
	
	\begin{minipage}{0.9\textwidth}
		\includegraphics[width=\textwidth]{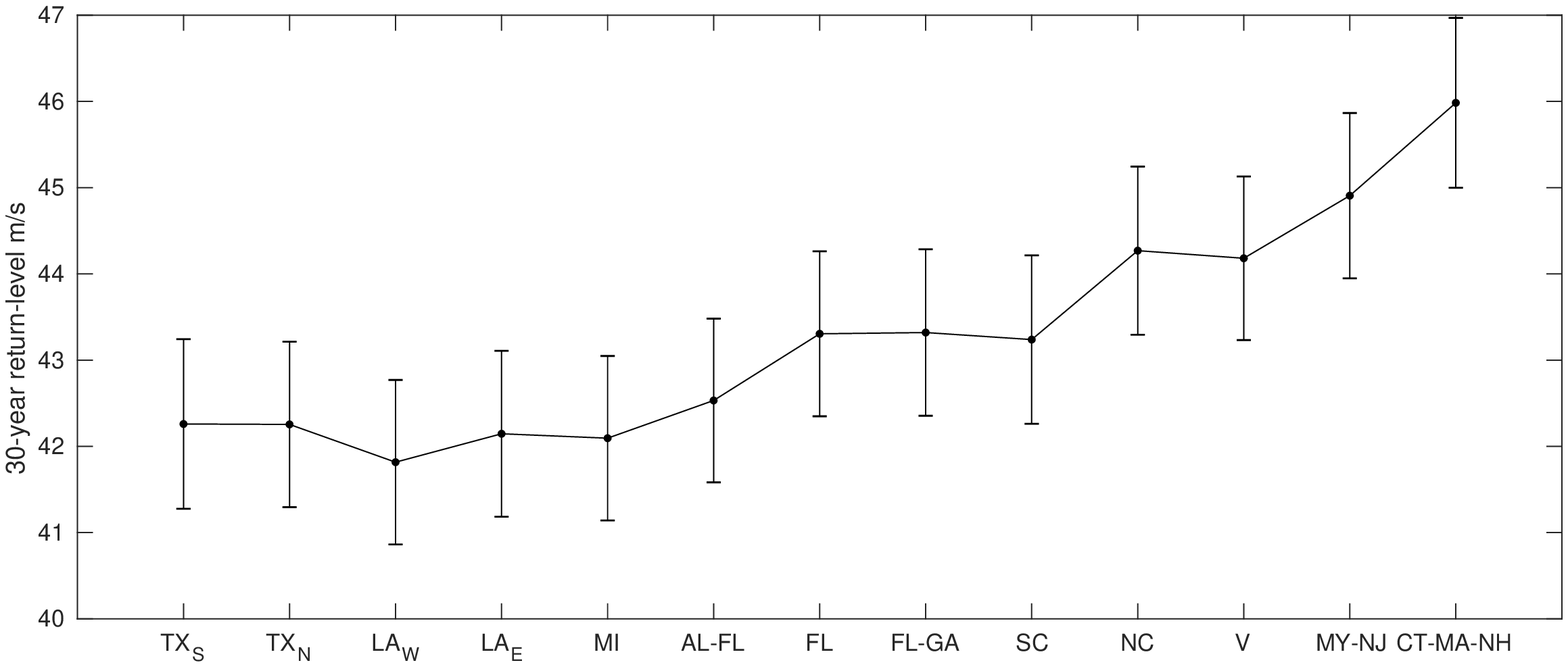}
	\end{minipage}
	
	\begin{minipage}{0.9\textwidth}
		\includegraphics[width=\textwidth]{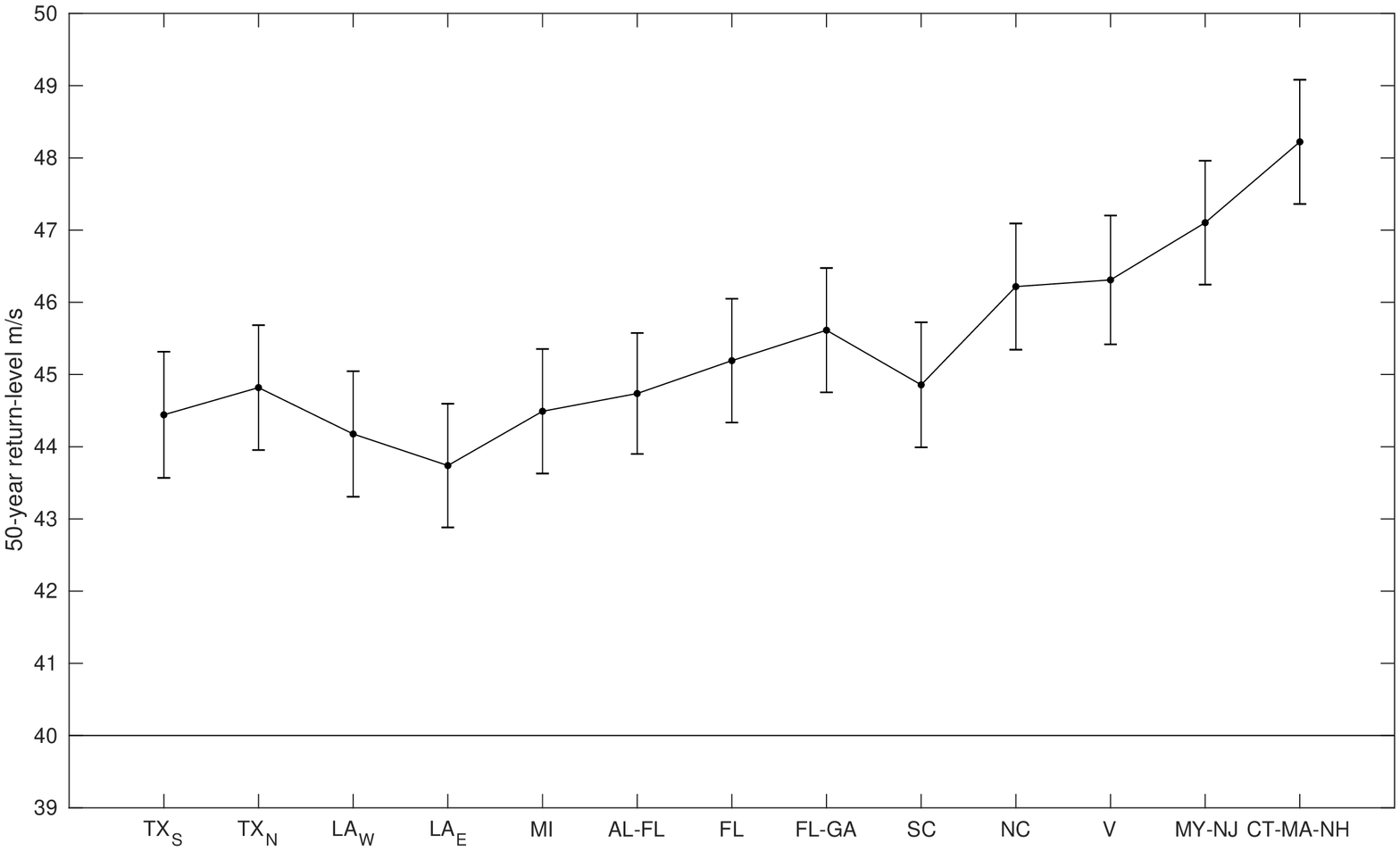}
	\end{minipage}
	\caption{Plots of the $20$-, $30$-, and $50$- year return-levels of maximum windspeed along the coast \textbf{estimated for year 2021}. Central estimates are the quantiles of the distribution of $6$-hourly windspeeds for $N = 1,000$ trials of $20$, $30$, and $50$ years of simulated hurricanes, respectively. Errorbars represent the $95\%$ confidence interval estimates from eq. \ref{95}. Solid horizontal line indicates the maximum windspeed return estimated from stationary models of previous literature. \label{ws_returns}}
\end{figure}

To estimate the confidence intervals of maximum windspeed return-levels, we independently sample from each of the parameter distributions to obtain 100 different combinations of parameters. We then run the simulation with each set of parameters for $2,000$ (e.g. $20$ years and $N=100$ trials), $3,000$ and $5,000$ years of hurricane simulations and estimate the $20$-, $30$- and $50$- year return-levels. Given that each of these simulations is independent, we are left with a sequence of quantile estimates (return-levels) coming from an i.i.d. sequence of maximum windspeeds for each coastal region. It show in \cite{K} that quantiles coming from i.i.d. sequences can be well-approximated by a normal distribution. Confidence intervals of each return-level are then estimated by assuming an underlying normal distribution so that,
\begin{equation}\label{95}
\text{CI}_{0.95} = \bigg[F^{-1}_{\text{norm}}(0.025)\frac{\sigma_{q}}{\sqrt{100}}, F^{-1}_{\text{norm}}(0.975)\frac{\sigma_{q}}{\sqrt{100}}\bigg]
\end{equation}
where $F^{-1}_{\text{norm}}$ is the inverse standard normal distribution and $\sigma_{q}$ is the estimated standard deviation of the 100 quantiles obtained from the 100 different parameter combinations. Quantiles to estimate $20$-, $30$-, and $50$-year return-levels of maximum windspeed for each coastal region and their 95\% estimated confidence intervals can be found in figure \ref{ws_returns}.

\section{Discussion}

The Wind Field Model \cite{NOAA} has provided a convenient way of calculating the maximum windspeed of a hurricane event at any given time along a track, provided the central pressure is known. According to this model, high maximum windspeeds are obtained for low central pressure measurements. It is shown in \cite{CC} that this relationship can be used as a guide for estimating returns of extremely high maximum windspeeds along the coast by appropriately modeling the pressure minima. They found using the HURDAT database that central pressure can be modeled by the generalized extreme value distribution with stationary location and scale parameters depending on the lifetime and latitude of the central pressure minima. The simulation results of \cite{CC} using a stationary model of central pressure minima are in good agreement with the other models and analyses of the decade \cite{B,H}. However, our investigation shows that this stationary model does not appropriately fit the central pressure minima of the updated HURDAT2 database. These poor fits can be almost entirely blamed on a time-dependent component of the scale and location parameters in the model. 

We have proposed a new, nonstationary model that accounts for this observed time-dependence in the location and scale parameter of the central pressure minima. Our model shows very reasonable fits to the true central pressure minima. Following a standard approach, we assume a Poisson distribution for yearly returns of hurricane events; however, we show that this model is also time-dependent with an exponentially increasing rate parameter for hurricane events with lifetimes greater than 6.25 days. We discuss this against current literature where stationarity of hurricane returns is assumed. We show that our models can reliably predict up to at least $50$-year returns for the central pressure minima without the need for updating by comparing the generated training set model against a test set and standardizing the central pressure minima using extreme value methods.

We have used our nonstationary model of central pressure minima and Poisson returns for yearly hurricanes in a more complex simulation to estimate $20$-, $30$-, and $50$- year return-levels of maximum windspeeds for sections along the US North Atlantic coastline. Compared to other analyses of maximum windspeed returns for coastal regions based on the HURDAT database like those in \cite{B,C,CC}, our model has two significant results: 1.) higher maximum windspeeds are expected for all regions along the US North Atlantic Coast and 2.) the highest maximum windspeeds occur in the northern part of the coast.

Specifically for landfalling hurricanes, we find a scale parameter for negative central pressure minima that is linearly increasing with time which suggests an expectation for higher-highs and lower-lows of central pressure minima. This phenomenon coupled with a general increase in the observed number of hurricane events can certainly lead to higher maximum windspeeds everywhere along the coast. 

An increase in the maximum windspeed for higher latitudes is actually a nontrivial observation because return-levels are affected by the number of hurricanes observed in a coastal region. In general, the number of observed hurricanes in the north tend to be lower. For example, it is well-known that the coastal region around Florida
has many more hurricane events than those regions along the north-eastern coast. In fact, we find this to be true in our simulations as well. So, one would expect to have a higher 20-year maximum wind-speed return-level for the Florida region than the northeast. On the other hand, translational velocity plays a critical role in the maximum wind-speed of a hurricane hitting the coastal region where translational velocity is always greater for higher latitudes.

\begin{remark}
The uneven latitudinal distribution of hurricane tracks is accounted for in the simulation by the use of historical tracks with some small added noise.
\end{remark}

\begin{figure}[h]
	\centering
	\includegraphics[width=0.8\textwidth]{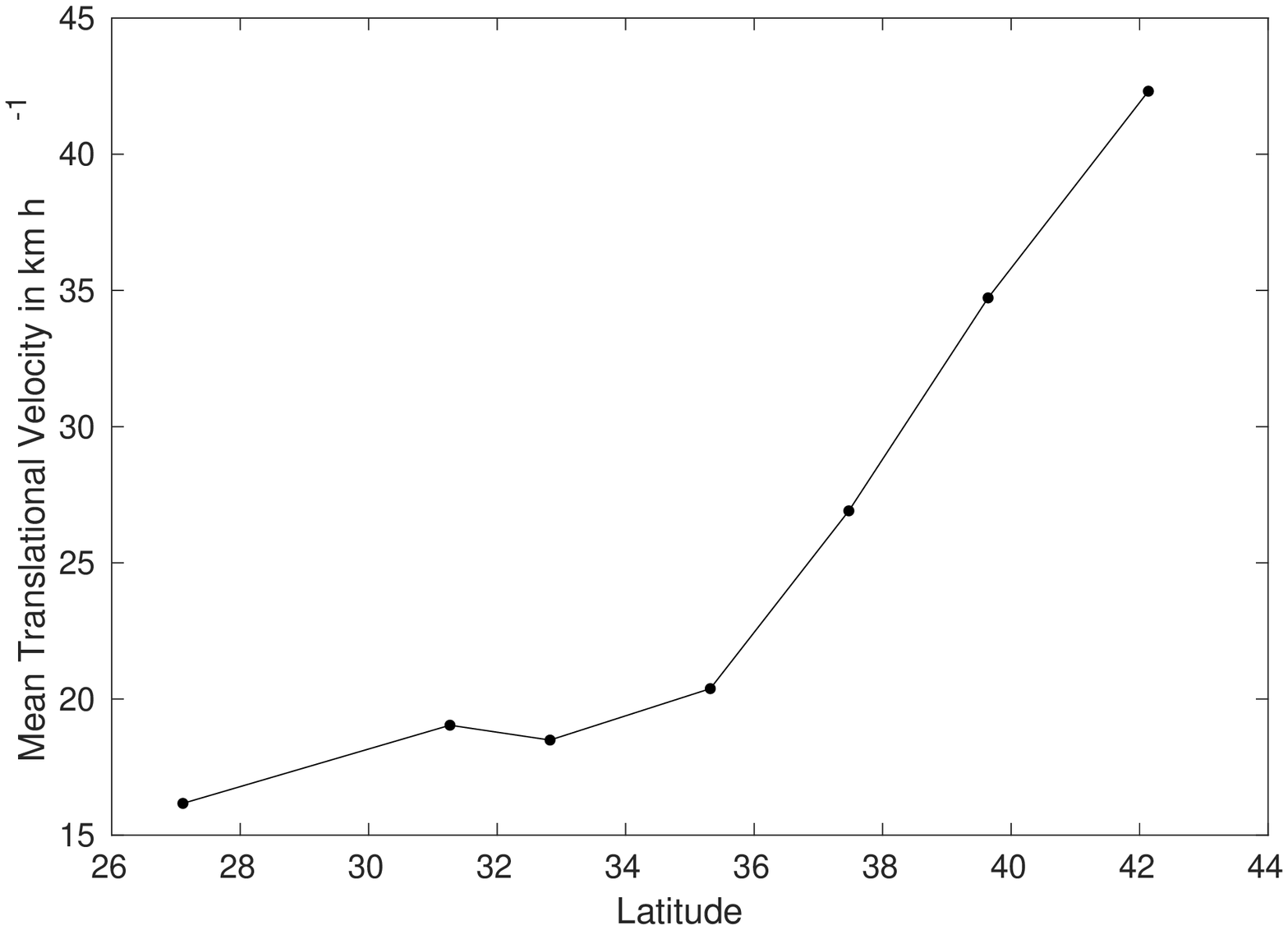}
	\caption{Average translational speed for simulated hurricanes from our model along the coast plotted against latitude. \label{tv}}
\end{figure}

We have tested our model to determine the cause of this northern increase in maximum windspeed and have found that translational velocity has the greatest influence over the observed trend. Furthermore, our simulated values of average translational velocity, estimated from the simulated hurricane tracks, almost identically follow those found in the literature \cite{Y}. We refer the reader to figure \ref{tv} for a plot of translational velocity over latitude. This result provides us with reasonable confidence in our model for coastal risk analysis.

\clearpage
\section*{Code and Data Availability}
The historical central pressure and track data are freely available on the NOAA website: https://www.nhc.noaa.gov/data/hurdat/hurdat2-1851-2020-052921.txt \cite{L2}. The code used to perform the nonstationary investigation, model fits, and wind-speed simulation was written in MATLAB. For access to this code, please contact the corresponding author.

\section*{Acknowledgements}
MC would like to thank MPIPKS for their hospitality where part of this work was completed. MN would like to thank the NSF for support on NSF-DMS Grants 1600780 and 2009923.

\clearpage
\appendix
\section{Numerical Simulation of Maximum Windspeeds for Hurricane Events}\label{simu}
\begin{figure}[h]
	\centering
	\includegraphics[scale=0.6]{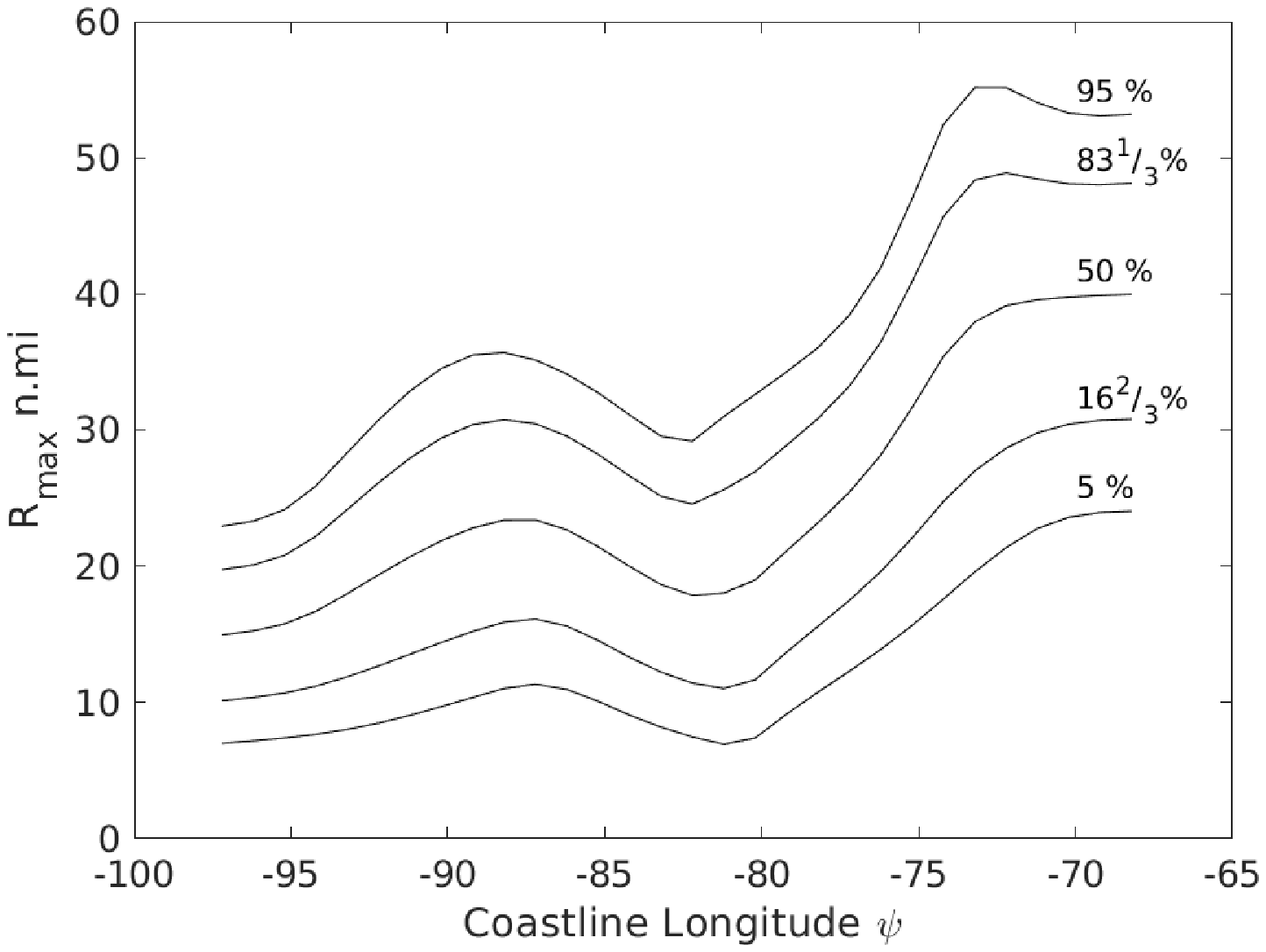}
		\caption{Observed quantiles of $R_{\max}$ as a function of the U.S. coastline longitude $\psi$ reproduced from fig. 37 and 38 of \cite{H}\label{rmax} and compared to fig. 7 of \cite{CC}. $R_{\max}$ increases as a function of the latitude as seen in the figure (increases and decreases along the U.S. coastline defined longitude). The model for $R_{\max}$ is created using the longitude $\psi$ while sampling is performed using the latitude $\phi$ because $R_{\max}$ is unique along $\phi$.}
\end{figure}

A summary of the scheme for simulating the timeseries of maximum windspeeds of a yearly sample of hurricane events is as follows. Adaptations from this investigation are marked with a *. Sampling tracks or timeseries refers to sampling from the 642 historical hurricane records from the HURDAT2 database. Densities and probabilities are estimated from the 642 historical hurricane records.
\begin{itemize}
	\item Sample the number of hurricane events to occur for a specified year from (\ref{poisson}) with rate parameter (\ref{rt}) *.
	\item For each hurricane event:
	\begin{itemize}
		\item Simulate a hurricane track by uniformly sampling a historical track and adding a small amount of spatial noise at each step sampled from $\mathcal{N}(0,\sigma^2)$ where $\sigma \simeq 100$ n.mi. 
		\item Uniformly sample a central pressure timeseries of a nonlandfalling hurricane.
		\item Time scale the central pressure timeseries to match the lifetime of the simulated
		track.
		\item For simulated landfalling hurricanes, run a Bernoulli trial with $p$ probability that the time of occurrence of the central pressure minimum occurring is at landfall, $t_{p_{\min}} = t_{\text{lf}}$, and $1-p$ otherwise. $p$ is estimated from the database. 
		\begin{itemize}
			\item If occurring before landfall, randomly sample the time of occurrence for central pressure minimum using the density of the ratio $t_{p_{\min}}/t_{\text{lf}}$ in the Appendix figure \ref{tratio}.
		\end{itemize}
		\item Randomly sample the central pressure minimum from the nonstationary model described by (\ref{ns_lf}) and (\ref{ns_nlf}) *.
		\item Randomly sample the central pressure range using the density model in \cite[Section 2.5]{CC} with fits in Table \ref{prange}.
		\item Add landfall effects described in \cite[Section 2.6]{CC} separately for inland and coastal hurricanes.
		\item Simulate the radius to maximum windspeeds $R_{\max}(\phi_t)$ by finding the distribution with coastal latitude $\phi_t$ illustrated in figure \ref{rmax}.
		\item Use the simulated central pressure timeseries $p_t$ and $R_{\max}(\phi_t)$ as inputs into the Wind Field Model described in (\ref{wfm}).
	\end{itemize}
\end{itemize}

\clearpage
\section{Additional Tables}

\begin{table}[ht]
	\centering
	\begin{tabular}{|llll|}
		\hline
		Type & $a$(se) & $b$(se) & $c$(se)\\
		\hline
		Landfalling & 872.33(25.48) & -0.87(0.03) & 11.77(0.32)\\
		Nonlandfalling & 829.50(19.77) & -0.82(0.02) & 7.47(0.19)\\
		\hline
	\end{tabular}
	\caption{Maximum likelihood estimates of the parameters in the normal distributional model from \cite[Section 2.5]{CC} for $p_{\text{range}}$.\label{prange}}
\end{table}

\clearpage
\section{Additional Figures}

\begin{figure}[h]
	\centering
	\includegraphics[scale=0.45]{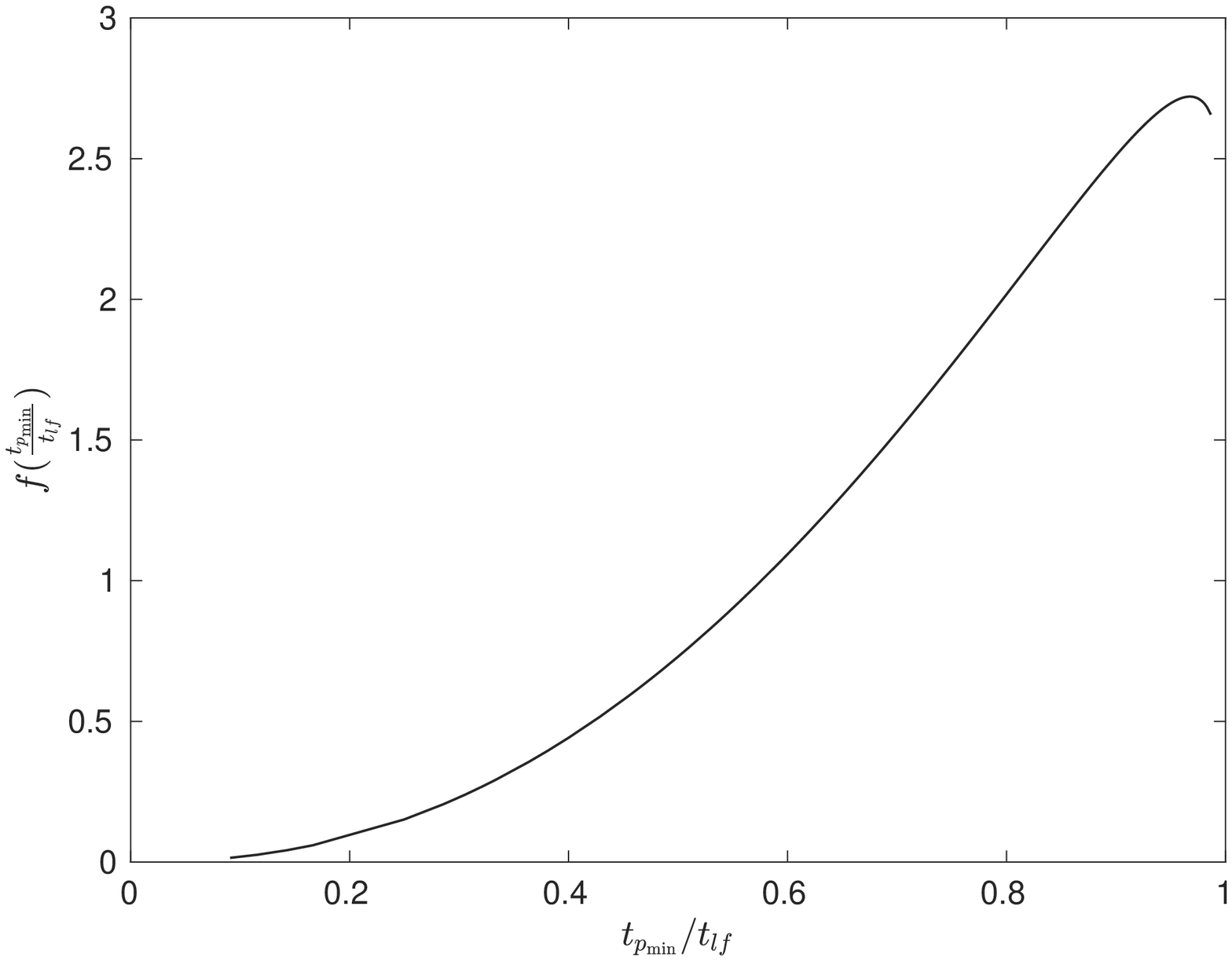}
	\caption{Empirical density function of the ratio $t_{\text{p}_{\text{min}}}/t_{\text{lf}}$ for landfalling hurricanes that make landfall after their central pressure minima occurs.\label{tratio}}
\end{figure}

\begin{figure}[h]
	\centering
	\begin{minipage}{\textwidth}
		\centering
		\includegraphics[scale=0.6]{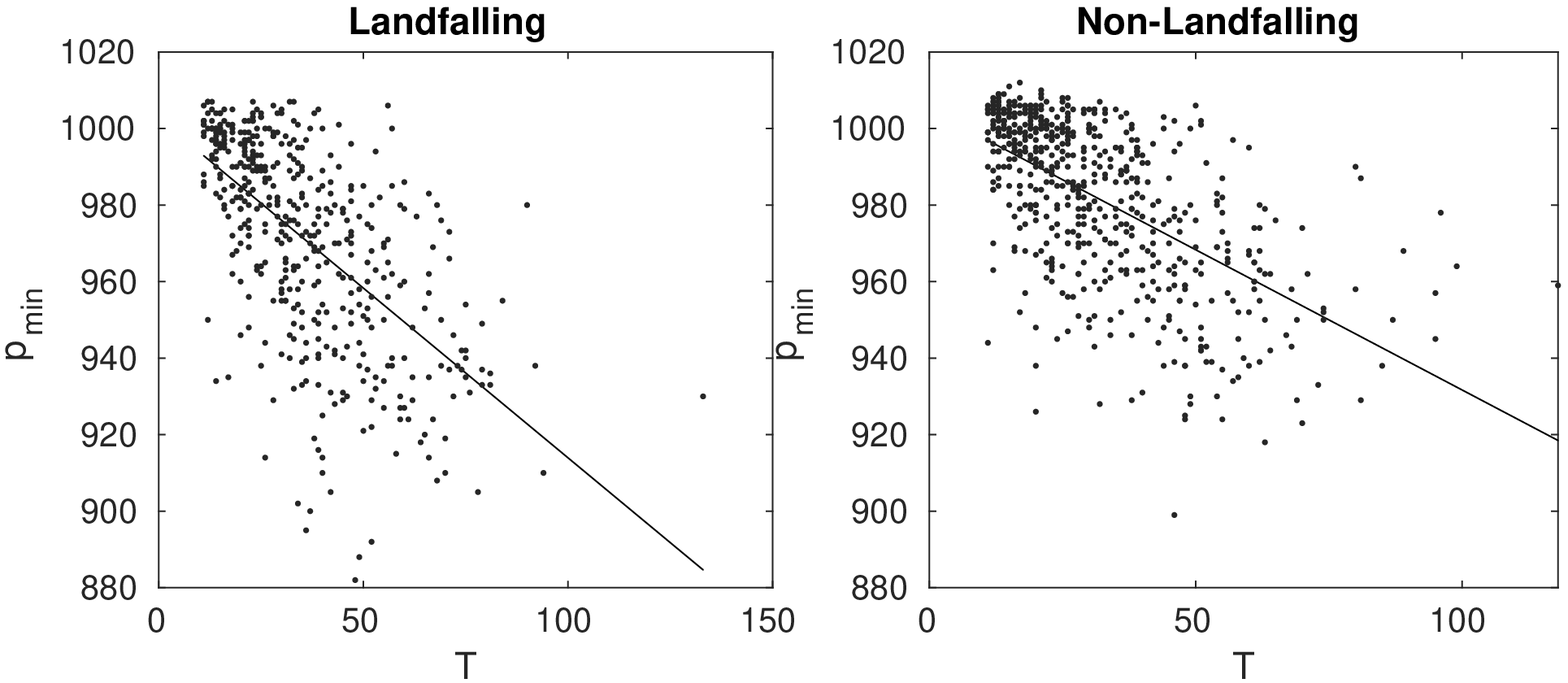}
		\caption*{(a)}
	\end{minipage}
	\begin{minipage}{\textwidth}
		\centering
		\includegraphics[scale=0.6]{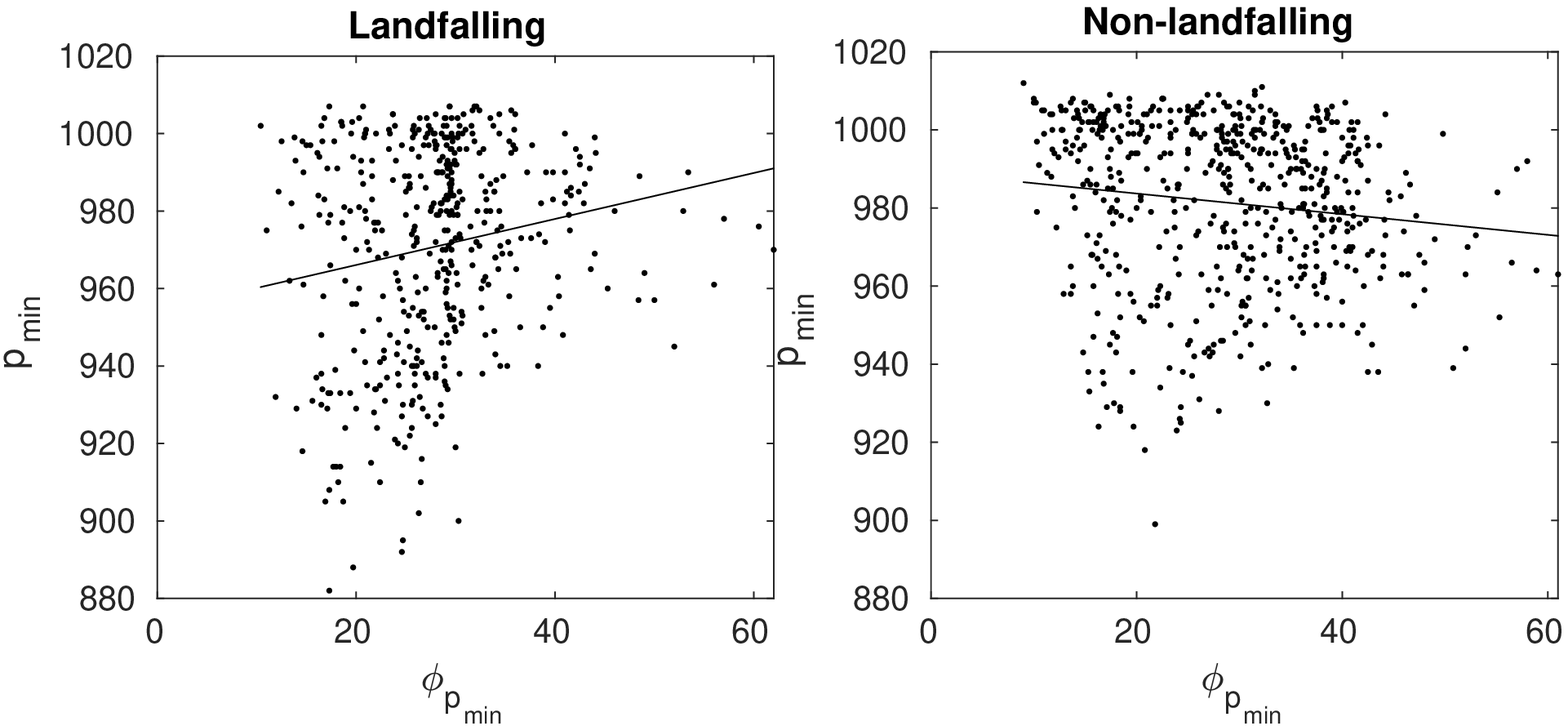}
		\caption*{(b)}
	\end{minipage}
	\caption{(a) Scatter plot of $p_{\min}$ against $T$ and (b) $\phi_{t_{\text{p}_{\min}}}$ for landfalling and nonlandfalling hurricanes with a line of best fit. Results indicate the the stationary model including dependence on $T$ and $\phi_{t_{p_{\min}}}$ is a reasonable starting point for forming the nonstationary model. \label{scatter}}
\end{figure}

\begin{figure}
	\centering
	\includegraphics[scale=0.6]{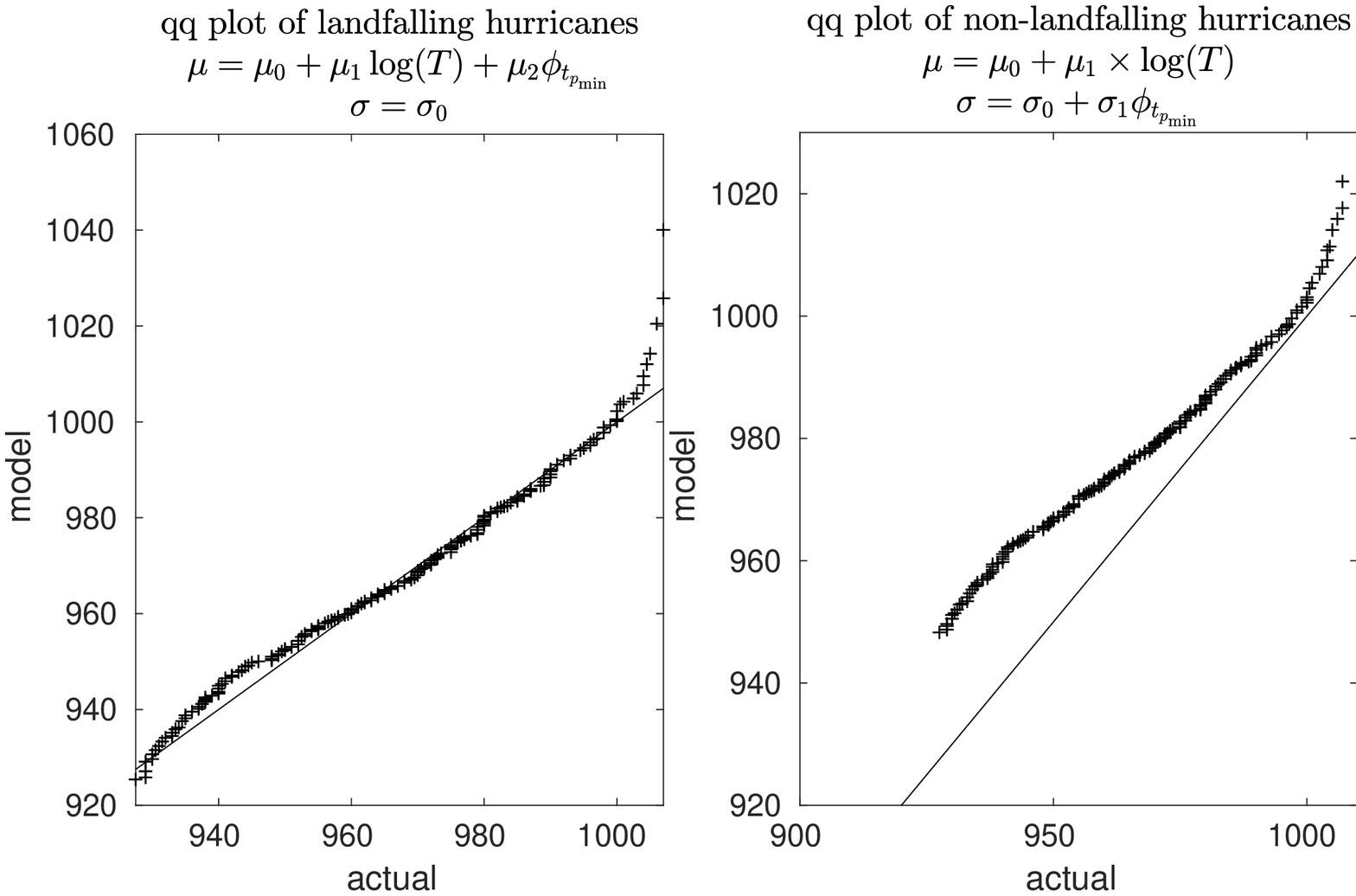}
	\caption{Quantile plots for landfalling and nonlandfalling hurricanes comparing the stationary model from \cite{CC} on the HURDAT2 database.\label{SFIT}}
\end{figure}

\begin{figure}
	\centering
	\includegraphics[scale=0.5]{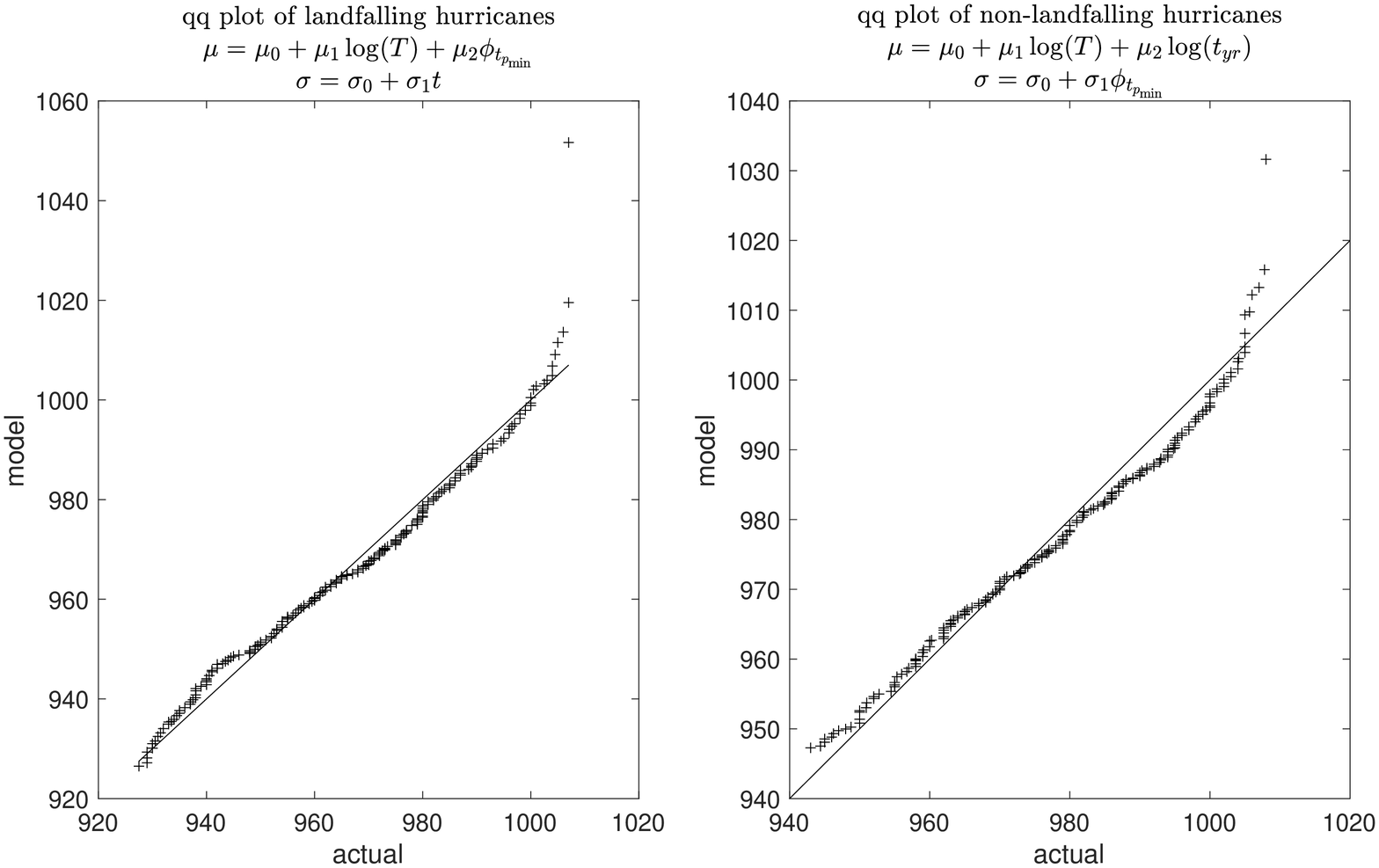}
	\caption{Quantile plots for landfalling and nonlandfalling hurricanes comparing the proposed nonstationary GEV model to data from the HURDAT2 database.\label{NSFIT}}
\end{figure}

\begin{figure}
	\centering
	\includegraphics[scale=0.12]{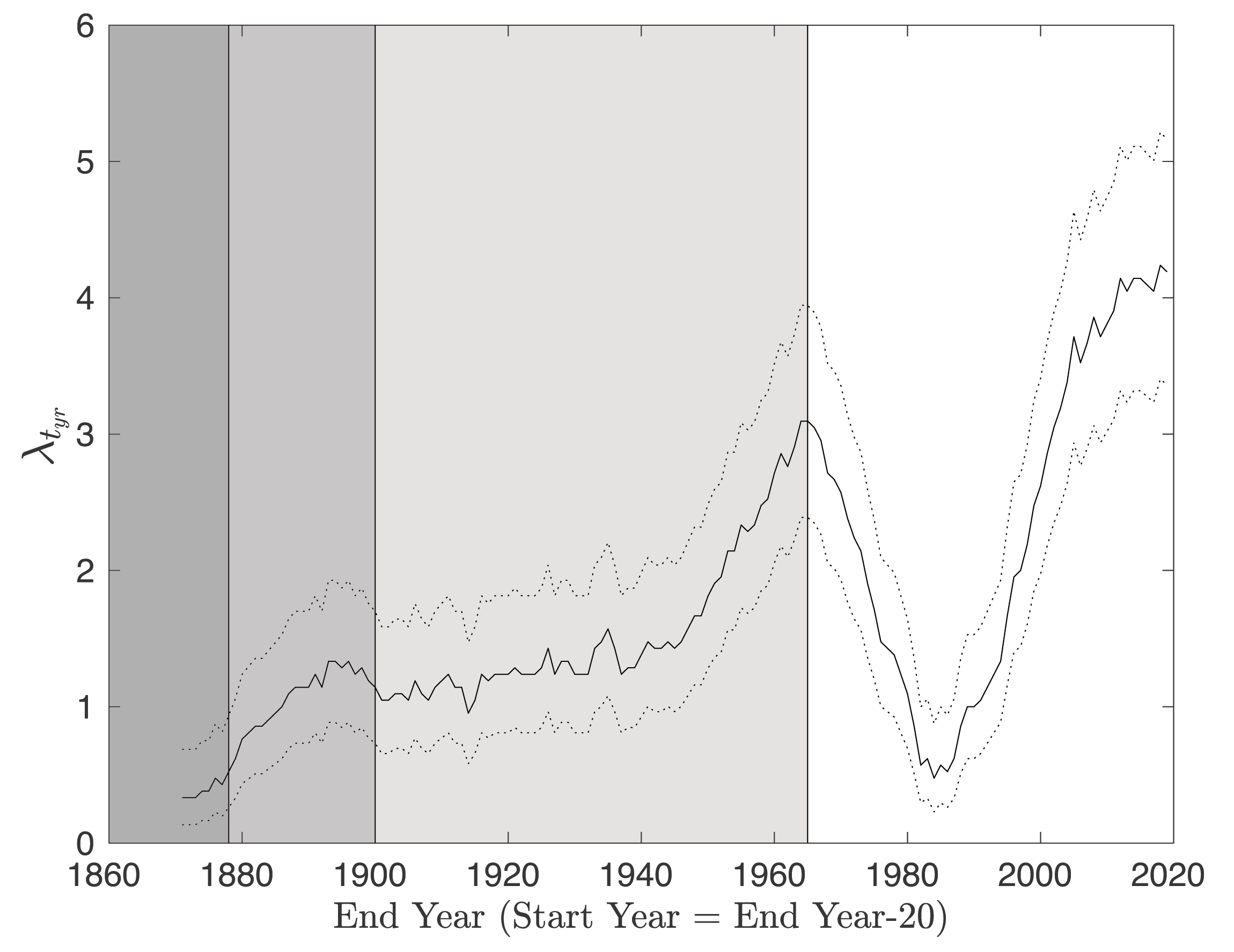}
	\caption{Likelihood estimate of the Poisson parameter for yearly landfalling hurricane event rates with lifetimes greater than 6.25 days. Estimates are taken over 20 year moving time windows. Standard errors are marked with dotted lines. Fitted exponential model is represented by a thick line. Grayed areas correspond to those in \cite{V3}: (1) 1878 - year when the U.S. Signal Corps began cataloging all Atlantic hurricanes (2) 1900 - year when the U.S. Coast was sufficiently well-populated for monitoring (3) modern-era with appropriate ship density. \label{POISSON_PARAM_LF}}
\end{figure}

\begin{figure}
	\centering
	\includegraphics[scale=0.14]{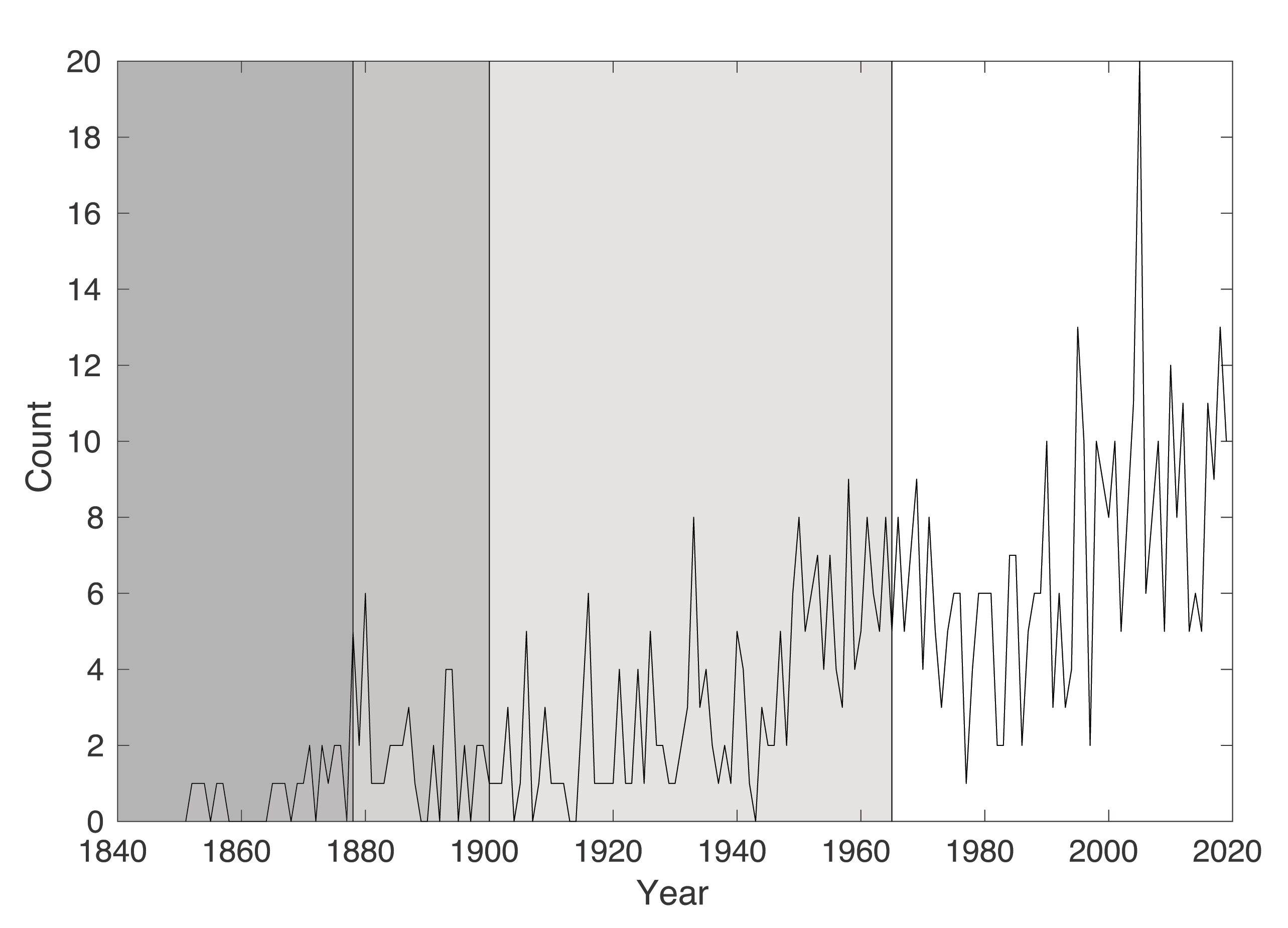}
	\caption{Number of total hurricane events of lifetimes greater than 6.25 days by year. Grayed areas correspond to those in \cite{V3}: (1) 1878 - year when the U.S. Signal Corps began cataloging all Atlantic hurricanes (2) 1900 - year when the U.S. Coast was sufficiently well-populated for monitoring (3) modern-era with appropriate ship density. \label{HURR_NUM}}
\end{figure}

\end{document}